\documentclass[12pt,a4paper,reqno]{amsart}
\usepackage{amssymb}

\usepackage{amscd}
\usepackage{enumerate}
\usepackage{graphicx}
\usepackage{color}
\numberwithin{equation}{section}

\usepackage[tableposition=top]{caption}
\usepackage{booktabs,dcolumn}

\theoremstyle{plain}

\newtheorem{theorem}{Theorem}[section]
\newtheorem{proposition}[theorem]{Proposition}

\newtheorem{lemma}[theorem]{Lemma}

\newtheorem{conjecture}[theorem]{Conjecture}

\newtheorem{remark}[theorem]{Remark}
\newtheorem*{remark*}{Remark}

\newtheorem{definition}[theorem]{Definition}

\newcommand\R{\mathbb{R}}
\newcommand\Z{\mathbb{Z}}
\newcommand\N{\mathbb{N}}
\newcommand\C{\mathbb{C}}
\newcommand\Q{\mathbb{Q}}
\newcommand\eps{\varepsilon}

\begin{document}

\title{An averaged form of Chowla's conjecture}

\author{Kaisa Matom\"aki}
\address{Department of Mathematics and Statistics \\
University of Turku, 20014 Turku\\
Finland}
\email{ksmato@utu.fi}

\author{Maksym Radziwi{\l}{\l}}
\address{Department of Mathematics\\
 Rutgers University, Hill Center for the Mathematical
Sciences \\
110 Frelinghuysen Rd., Piscataway, NJ 08854-8019}
\email{maksym.radziwill@gmail.com}

\author{Terence Tao}
\address{Department of Mathematics, UCLA\\
405 Hilgard Ave\\
Los Angeles CA 90095\\
USA}
\email{tao@math.ucla.edu}

\begin{abstract}  Let $\lambda$ denote the Liouville function.  A well known conjecture of Chowla asserts that for any distinct natural numbers $h_1,\dots,h_k$, one has $\sum_{1 \leq n \leq X} \lambda(n+h_1) \dotsm \lambda(n+h_k) = o(X)$ as $X \to \infty$.  This conjecture remains unproven for any $h_1,\dots,h_k$ with $k \geq 2$.  In this paper, using the recent results of the first two authors on mean values of multiplicative functions in short intervals, combined with an argument of Katai and Bourgain-Sarnak-Ziegler, we establish an averaged version of this conjecture, namely 
$$ \sum_{h_1,\dots,h_k \leq H} \left|\sum_{1 \leq n \leq X} \lambda(n+h_1) \dotsm \lambda(n+h_k)\right| = o(H^kX)$$
as $X \to \infty$ whenever $H = H(X) \leq X$ goes to infinity as $X \to \infty$, and $k$ is fixed.  Related to this, we give the exponential sum estimate
$$ \int_0^X \left|\sum_{x \leq n \leq x+H} \lambda(n) e(\alpha n)\right| dx = o( HX )$$
as $X \to \infty$ uniformly for all $\alpha \in \R$, with $H$ as before.  Our arguments in fact give quantitative bounds on the decay rate (roughly on the order of $\frac{\log\log H}{\log H}$), and extend to more general bounded multiplicative functions than the Liouville function, yielding an averaged form of a (corrected) conjecture of Elliott.
\end{abstract}

\maketitle


\section{Introduction}

Let $\lambda: \N \to \{-1,+1\}$ be the Liouville function, that is to say the completely multiplicative function such that $\lambda(p)=-1$ for all primes $p$.  The prime number theorem implies that \footnote{See Section \ref{notation-sec} below for our asymptotic notation conventions.}
$$ \sum_{1 \leq n \leq X} \lambda(n) = o(X)$$
as $X \to \infty$.  More generally, a famous conjecture of Chowla \cite{Chowla-book} asserts that for any distinct natural numbers $h_1,\dots,h_k$, one has
\begin{equation}\label{chowla-conj}
 \sum_{1 \leq n \leq X} \lambda(n+h_1) \dotsm \lambda(n+h_k) = o(X)
\end{equation}
as $X \to \infty$.  

Chowla's conjecture remains open for any $h_1,\dots,h_k$ with $k \geq 2$. 
 Our first main theorem establishes an averaged form of this conjecture:

\begin{theorem}[Chowla's conjecture on average]\label{main}  For any natural number $k$, and any $10 \leq H \leq X$, we have
\begin{equation}\label{k-cor}
 \sum_{1 \leq h_1,\dots,h_k \leq H} \left|\sum_{1 \leq n \leq X} \lambda(n+h_1) \dotsm \lambda(n+h_k)\right| \ll k \left( \frac{\log\log H}{\log H} + \frac{1}{\log^{1/3000} X} \right) H^k X.
\end{equation}
In fact, we have the slightly stronger bound
\begin{equation}\label{k-cor-2}
 \sum_{1 \leq h_2,\dots,h_k \leq H} \left|\sum_{1 \leq n \leq X} \lambda(n) \lambda(n+h_2) \dotsm \lambda(n+h_k)\right| \ll k \left( \frac{\log\log H}{\log H} + \frac{1}{\log^{1/3000} X} \right) H^{k-1} X.
\end{equation}
\end{theorem}
In the case $k = 2$ our result implies that
$$
\sum_{1 \leq h \leq H} \Big | \sum_{1 \leq n \leq X} \lambda (n) \lambda(n + h) \Big | = o(HX)
$$
provided that $H \rightarrow \infty$ arbitrarily slowly with $X \rightarrow \infty$ (and $H \leq X$). Note that the $k=2$ case of Chowla's conjecture is equivalent to the above asymptotic holding in the case that $H$ is bounded rather than going to infinity.

In fact, we have a more precise bound than \eqref{k-cor} (or \eqref{k-cor-2}) that gives more control on the exceptional tuples $(h_1,\dots,h_k)$ for which the sums $\sum_{1 \leq n \leq X} \lambda(n+h_1) \dotsm \lambda(n+h_k)$ are large; see Remark \ref{improv} below. In particular in the special case $k=2$ we get the following result.
\begin{theorem}
\label{th:quant}
Let $\delta \in (0, 1]$ be fixed. There is a large but fixed $H = H(\delta)$ such that, for all large enough $X$,
\begin{equation}
\label{eq:quanteq}
\left|\sum_{1 \leq n \leq X}  \lambda(n) \lambda(n+h) \right| \leq \delta X
\end{equation}
for all but at most $H^{1-\frac{\delta}{5000}}$ integers $|h| \leq H$.
\end{theorem}
One can also replace the ranges $1 \leq h_j \leq H$ in Theorem~\ref{main} by $b_j+1 \leq h_j \leq b_j+H$ for any $b_j=O(X)$; see Theorem \ref{elliott-thm} below.

The exponents $1/3000$ and $1/5000$ in the above theorems may certainly be improved, but we did not attempt to optimize the constants here.  However, our methods cannot produce a gain much larger than $\frac{1}{\log H}$, as one would then have to somehow control $\lambda$ on numbers that are not divisible by any prime less than $H$, at which point we are no longer able to exploit the averaging in the $h_1,\dots,h_k$ parameters.  It would be of particular interest to obtain a gain of more than $\frac{1}{\log X}$, as one could then potentially localize $\lambda$ to primes and obtain some version of the prime tuples conjecture when the $h_1,\dots,h_k$ parameters are averaged over short intervals, but this is well beyond the capability of our methods.  (If instead one is allowed to average the $h_1,\dots,h_k$ over long intervals (of scale comparable to $X$), one can obtain various averaged forms of the prime tuples conjecture and its relatives, by rather different methods to those used here; see \cite{balog}, \cite{mikawa}, \cite{kawada}, \cite{kawada2}, \cite{green-tao}.)

Theorem \ref{main} is closely related to the following averaged short exponential sum estimate, which may be of independent interest.

\begin{theorem}[Exponential sum estimate]\label{main-2}  For any $10 \leq H \leq X$, one has
$$ \sup_{\alpha \in \R} \int_0^{X} \left| \sum_{x \leq n \leq x+H} \lambda(n) e(\alpha n) \right|\ dx \ll \left( \frac{\log\log H}{\log H} + \frac{1}{\log^{1/700} X} \right) HX.$$
\end{theorem}

Actually, for technical reasons it is convenient to prove a sharper version of Theorem \ref{main-2} in which the Liouville function has been restricted to those numbers that have ``typical'' factorization; see Theorem \ref{second}.  This sharper version will then be used to establish Theorem \ref{main}.   

The relationship between Theorem \ref{main} and Theorem \ref{main-2} stems from the following Fourier-analytic identity:

\begin{lemma}[Fourier identity]\label{fourier}  If $f: \Z \to \C$ is a function supported on a finite set, and $H > 0$, then
$$ \int_{{\mathbb T}} \left(\int_\R \left|\sum_{x \leq n \leq x+H} f(n) e(\alpha n)\right|^2\ dx\right)^2\ d\alpha = \sum_{|h| \leq H} (H-|h|)^2 \left|\sum_n f(n) \overline{f}(n+h)\right|^2.$$
\end{lemma}

\begin{proof}  Using the Fourier identity $\int_{{\mathbb T}} e( n \alpha )\ d\alpha = 1_{n=0}$, we can expand the left-hand side as
$$ \sum_{n, n', m, m'} f(n) \overline{f}(n') f(m) \overline{f}(m') 1_{n+m-n'-m'=0} \left(\int_\R 1_{x \leq n,n' \leq x+H}\ dx\right) \left(\int_\R 1_{y \leq m,m' \leq y+H}\ dy\right).$$
Writing $n' = n+h$, we see that both integrals are equal to $H-|h|$ if $|h| \leq H$, and vanish otherwise.  The claim follows.
\end{proof}

Theorem \ref{main-2} may be compared with the classical estimate
$$ \sup_{\alpha \in \R} \left|\sum_{1 \leq n \leq X} \lambda(n) e(\alpha n)\right| \ll_A X \log^{-A} X $$
of Davenport \cite{davenport}, valid for any $A > 0$.  Indeed, one can view Theorem \ref{main-2} as asserting that a weak form of Davenport's estimate holds on average in short intervals.  It would be of interest to also obtain non-trivial bounds on the larger quantity
\begin{equation}\label{larger}
 \int_0^{X} \sup_{\alpha \in \R} \left| \sum_{x \leq n \leq x+H} \lambda(n) e(\alpha n) \right|\ dx 
\end{equation}
but this appears difficult to establish with our methods.

As with other applications of the circle method, our proof of Theorem \ref{main-2} splits into two cases, depending on whether the quantity $\alpha$ is on ``major arc'' or on ``minor arc''.  In the ``major arc'' case we are able to use the recent results of the first two authors \cite{MR} on the average size of mean values of multiplicative functions on short intervals.  Actually, in order to handle the presence of complex Dirichlet characters, we need to extend the results in \cite{MR} to complex-valued multiplicative functions rather than real-valued ones; this is accomplished in an appendix to this paper (Appendix \ref{complex-app}).  In the ``minor arc'' case we use a variant of the arguments of Katai \cite{Katai} and Bourgain-Sarnak-Ziegler \cite{BSZ} (see also the earlier works of Montgomery-Vaughan \cite{MontgomeryVaughan} and 
Daboussi-Delange \cite{DaboussiDelange}) to obtain the required cancellation.  One innovation here is to rely on a combinatorial identity of Ramar\'e (also used in \cite{MR}) as a substitute for the Turan-Kubilius inequality, as this leads to superior quantitative estimates (particularly if one first restricts the variable $n$ to have a ``typical'' prime factorization).

\subsection{Extension to more general multiplicative functions}

Define a \emph{$1$-bounded multiplicative function} to be a multiplicative function $f: \N \to \C$ such that $|f(n)| \leq 1$ for all $n \in \N$.  Given two $1$-bounded multiplicative functions $f, g$ and a parameter $X \geq 1$, we define the distance ${\mathbb D}(f,g;X) \in [0,+\infty)$ by the formula
$$ {\mathbb D}(f,g;X) := \left(\sum_{p \leq X} \frac{1 - \operatorname{Re}( f(p) \overline{g(p)} )}{p}\right)^{1/2}.$$
This is known to give a (pseudo-)metric on $1$-bounded multiplicative functions; see \cite[Lemma 3.1]{GS}.  We also define the asymptotic counterpart ${\mathbb D}(f,g;\infty) \in [0,+\infty]$ by the formula
$$ {\mathbb D}(f,g;\infty) := \left(\sum_p \frac{1 - \operatorname{Re}( f(p) \overline{g(p)} )}{p}\right)^{1/2}.$$
We informally say that $f$ \emph{pretends} to be $g$ if ${\mathbb D}(f,g;X)$ (or ${\mathbb D}(f,g;\infty)$) is small (or finite).  

For any $1$-bounded multiplicative function $g$ and real number $X>1$, we introduce the quantity
\begin{equation} 
\label{eq:MgXdef}
M(g;X) := \inf_{|t| \leq X}  {\mathbb D}( g, n \mapsto n^{it}; X )^2,
\end{equation}
and then the more general quantity
\begin{align*}
M(g; X, Q) &:= \inf_{q \leq Q; \chi\ (q)} M(g \overline{\chi}; X) \\
&= \inf_{|t| \leq X; q \leq Q; \chi\ (q)} {\mathbb D}( g, n \mapsto \chi(n) n^{it}; X )^2,
\end{align*}
where $\chi$ ranges over all Dirichlet characters of modulus $q \leq Q$.  Informally, $M(g;X)$ is small when $g$ pretends to be like a multiplicative character $n \mapsto n^{it}$, and $M(g; X,Q)$ is small when $g$ pretends to be like a twisted Dirichlet character of modulus at most $Q$ and twist of height at most $X$.  We also define the asymptotic counterpart
$$ M(g;\infty,\infty) = \inf_{\chi, t} {\mathbb D}( g, n \mapsto \chi(n) n^{it}; \infty )^2$$
where $\chi$ now ranges over all Dirichlet characters and $t$ ranges over all real numbers.

In \cite[Conjecture II]{elliott}, Elliott proposed the following more general form of Chowla's conjecture, which we phrase here in contrapositive form.

\begin{conjecture}[Elliott's conjecture]\label{elliott-conj}  Let $g_1,\dots,g_k: \N \to \C$ be $1$-bounded multiplicative functions, and let $a_1,\dots,a_k,b_1,\dots,b_k$ be natural numbers such that any two of the $(a_1,b_1),\dots,(a_k,b_k)$ are linearly independent in $\Q^2$.  Suppose that there is an index $1 \leq j_0 \leq k$ such that
\begin{equation}\label{mgi}
 M( g_{j_0}; \infty, \infty ) = \infty.
\end{equation}
Then
\begin{equation}\label{el}
 \sum_{1 \leq n \leq X} \prod_{j=1}^k g_j(a_j n + b_j) = o(X)
\end{equation}
as $X \to \infty$.
\end{conjecture}

Informally, this conjecture asserts that for pairwise linearly independent $(a_1,b_1),\dots,(a_k,b_k)$ and any $1$-bounded multiplicative $g_1,\dots,g_k$, one has the asymptotic \eqref{el} as $X \to \infty$, unless each of the $g_j$ pretends to be a twisted Dirichlet character $n \mapsto \chi_j(n) n^{it_j}$.  Note that some condition of this form is necessary, since if $g(n)$ is equal to $\chi(n) n^{it}$ then $g(n) \overline{g(n+h)}$ will be biased to be positive for large $n$, if $h$ is fixed and divisible by the modulus $q$ of $\chi$; one also expects some bias when $h$ is not divisible by this modulus since the sums $\sum_{n \in \Z/q\Z} \chi(n) \overline{\chi(n+h)}$ do not vanish in general.  From the prime number theorem in arithmetic progressions it follows that
$$
M( \lambda; \infty, \infty) = \infty, 
$$
so Elliott's conjecture implies Chowla's conjecture \eqref{chowla-conj}.

When one allows the functions $g_j$ to be complex-valued rather than real-valued, Elliott's conjecture turns out to be false on a technicality; one can choose $1$-bounded multiplicative functions $g_j$ which are arbitrarily close at various scales to a sequence of functions of the form $n \mapsto n^{it_m}$ (which allows one to violate \eqref{el}) without \emph{globally} pretending to be $n^{it}$ (or $\chi(n) n^{it}$) for any \emph{fixed} $t$; we present this counterexample in Appendix \ref{counter}.  However, this counterexample can be removed by replacing \eqref{mgi} with the stronger condition that
\begin{equation}\label{geo}
 M(g_{j_0}; X, Q) \to \infty
\end{equation}
as $X \to \infty$ for each fixed $Q$.  In the real-valued case, \eqref{geo} and \eqref{mgi} are equivalent by a triangle inequality argument of Granville and Soundararajan which we give in Appendix \ref{gs-sec}.

As evidence for the corrected form of Conjecture \ref{elliott-conj} (in both the real-valued and complex-valued cases), we present the following averaged form of that conjecture:

\begin{theorem}[Elliott's conjecture on average]\label{elliott-thm}  Let $10 \leq H \leq X$ and $A \geq 1$.
Let $g_1,\dots,g_k: \N \to \C$ be $1$-bounded functions, and let $a_1,\dots,a_k,b_1,\dots,b_k$ be natural numbers with $a_j \leq A$ and $b_j \leq AX$ for $j=1,\dots,k$.  Let $1 \leq j_0 \leq k$, and suppose that $g_{j_0}$ is multiplicative.  Then one has
\begin{align}\label{ag}
 \sum_{1 \leq h_1,\dots,h_k \leq H} & \left|\sum_{1 \leq n \leq X} \prod_{j=1}^k  g_j(a_j n + b_j + h_j)\right| \\ \nonumber & \ll A^2 k
\left( \exp( - M / 80) + \frac{\log\log H}{\log H} + \frac{1}{\log^{1/3000} X} \right) H^k X
\end{align}
where
$$ M := M(g_{j_0};10AX, Q)$$
and 
$$ Q := \min( \log^{1/125} X, \log^{20} H ).$$
In fact, we have the slightly stronger bound
\begin{align}\label{bg}
 \sum_{1 \leq h_2,\dots,h_k \leq H} & \left|\sum_{1 \leq n \leq X} g_1(a_1 n + b_1) \prod_{j=2}^k g_j(a_j n + b_j + h_j)\right| \\ & \nonumber \ll A^2
k \left( \exp( - M / 80) + \frac{\log\log H}{\log H} + \frac{1}{\log^{1/3000} X} \right) H^{k-1} X.
\end{align}
\end{theorem}

Note that if $a_1,\dots,a_k,b_1,\dots,b_k$ are fixed, $g_{j_0}$ is independent of $X$ and obeys the condition \eqref{geo} for any fixed $Q$, and $H = H(X)$ is chosen to go to infinity arbitrarily slowly as $X \to \infty$, then the quantity $M$ in the above theorem goes to infinity (note that $M(g;X,Q)$ is non-decreasing in $Q$), and \eqref{bg} then implies an averaged form of the asymptotic \eqref{el}.  Thus Theorem \ref{elliott-thm} is indeed an averaged form of the corrected form of Conjecture \ref{elliott-conj}.  (We discovered the counterexample in Appendix \ref{counter} while trying to interpret Theorem \ref{elliott-thm} as an averaged version of the original form of Conjecture \ref{elliott-conj}.)  Interestingly, only one of the functions $g_1,\dots,g_k$ in Theorem \ref{elliott-thm} is required to be multiplicative\footnote{We thank the referee for observing this fact.  In a previous version of this paper, all of the $g_j$ were required to be multiplicative.}; one can use a van der Corput argument to reduce matters to obtaining cancellation for a sum roughly of the form $\sum_{h \leq H} |\sum_{1 \leq n \leq X} g_{j_0}(n) \overline{g_{j_0}(n+h)}|^2$, which can then be treated using Lemma \ref{fourier}.

For $g(n) = \lambda(n)$ and $X,Q,M$ as in the above theorem, one obtains, for every $\varepsilon > 0$, the bound
\begin{equation}
\label{eq:MestLiou}
M \geq \inf_{|t| \leq X; q \leq Q; \chi\ (q)} \sum_{\exp((\log X)^{2/3+\varepsilon}) \leq p \leq X} \frac{1+ \operatorname{Re} \chi(p) p^{it}}{p} \geq \left(\frac{1}{3}-\varepsilon\right) \log \log X + O(1)
\end{equation}
where the last inequality is established via standard methods from the Vinogradov-Korobov type zero-free region 
\[
\left\{ \sigma+it: \sigma > 1-\frac{c}{\max\{\log q, (\log (3+|t|))^{2/3} (\log \log (3+|t|))^{1/3}\}} \right\}
\] 
for $L(s, \chi)$ and some absolute constant $c > 0$, which applies since $\chi$ has conductor $q \leq (\log X)^{1/125}$ (so that there are no exceptional zeros), see \cite[\S 9.5]{Montgomery}. Hence Theorem~\ref{elliott-thm} implies Theorem~\ref{main}.  The same argument gives Theorem \ref{main} when the Liouville function $\lambda$ is replaced by the M\"obius function $\mu$.  We remark that as our arguments make no use of exceptional zeroes, all the implied constants in our theorems are effective.

We also have a generalized form of Theorem \ref{main-2}:

\begin{theorem}[Exponential sum estimate]\label{exp-est}  Let $X \geq H \geq 10$ and let $g$ be a $1$-bounded multiplicative function.
Then 
\begin{align*}
\sup_{\alpha \in {\mathbb T}} \int_0^{X} \Big |\sum_{x \leq n \leq x+H} g(n) & e(\alpha n) \Big |\ dx    \\ & \ll 
\left( \exp( -M(g; X, Q) / 20 ) + \frac{\log \log H}{\log H} + \frac{1}{\log^{1/700} X} \right) HX
\end{align*}
where
$$ Q := \min( \log^{1/125} X, \log^{5} H ).$$
\end{theorem}

By~\eqref{eq:MestLiou}, Theorem \ref{exp-est} implies Theorem \ref{main-2}.  

\begin{remark}  In the recent paper \cite{FH}, a different averaged form of Elliott's conjecture is established, in which one uses fewer averaging parameters $h_i$ than in Theorem \ref{elliott-thm} (indeed, one can average over just a single such parameter, provided that the linear parts of the forms are independent), but the averaging parameters range over a long range (comparable to $X$) rather than on the short range given here.  The methods of proof are rather different (in particular, the arguments in \cite{FH} rely on higher order Fourier analysis).  In the long-range averaged situation considered in \cite{FH}, the counterexample in Appendix \ref{counter} does not apply, and one can use the original form of Elliott's conjecture in place of the corrected version.  It may be possible to combine the results here with those in \cite{FH} to obtain an averaged version of the Chowla or Elliott's conjecture in which the number of averaging parameters is small, and the averaging is over a short range, but this seems to require non-trivial estimates on quantities such as \eqref{larger}, which we are currently unable to handle.
\end{remark}

\begin{remark}  Theorem \ref{elliott-thm} suggests that in order to make the correlation
$$ \sum_{1 \leq n \leq X} \prod_{j=1}^k  g_j(a_j n + b_j )$$
significantly smaller than $X$, one should have $M(g_{j_0};AX,A)$ large for some moderately large $A$ and some $1 \leq j_0 \leq k$.  This appears to be the right condition when $k=2$, but for larger values of $k$ it appears that one in fact should require that $ M(g_{j_0};AX^{k-1},A)$ is large; that is to say, one may conjecture the bound
$$ |\sum_{1 \leq n \leq X} \prod_{j=1}^k  g_j(a_j n + b_j )| \leq \eps X$$
whenever $M(g_{j_0};AX^{k-1},A) \geq A$ for some $A$ sufficiently large depending on $\eps$, $k$, and the $a_j,b_j$, assuming that $X$ is sufficiently large depending on $A,\eps,k,a_j,b_j$, that one has $a_i b_j - a_j b_i \neq 0$ for all $i \neq j$, and the $g_j$ are all $1$-bounded multiplicative functions.  For instance, consider the $k=3$ correlation
$$ \sum_{1 \leq n \leq X} g_1(n) g_2(n+1) g_3(n+2).$$
This sum is large in the case $g_1(n) := n^{it}$, $g_2(n) := n^{-2it}$, $g_3(n) := n^{it}$ for $t = o(X^2)$, as can be seen by applying Taylor expansion to second order to the function $x \mapsto t \log x$ around $x = n$.  If $t$ is much larger than $X$, then the quantities $M(g_{j_0};AX,A)$ are large for $j_0=1,2,3$, but $M(g_{j_0};AX^2,A)$ are small, and so needs a lower bound on $M(g_{j_0};AX^2,A)$ rather than $M(g_{j_0};AX,A)$ to ensure the smallness of this correlation.  A similar example can be constructed for higher values of $k$.  A related computation also shows that if one wishes to move the supremum in $\alpha$ in Theorem \ref{exp-est} inside the integration in $x$ (as in \eqref{larger}), one will need a lower bound on $M(g; X^2, Q)$ rather than just $M(g; X, Q)$.
\end{remark}

\subsection{Acknowledgments}

TT was supported by a Simons Investigator grant, the James and Carol Collins Chair, the Mathematical Analysis \&
Application Research Fund Endowment, and by NSF grant DMS-1266164. The authors thank Andrew Granville and the anonymous referee for useful comments and corrections.
We thank Fei Wei for pointing out the issue with the proof of Proposition A.3 in the published version of the paper. We thank Alisa Sedunova and Ke Wang for pointing out the correction to Ramar{\'e}'s identity.

\subsection{Notation}\label{notation-sec}

Our asymptotic notation conventions are as follows.  We use $X \ll Y$, $Y \gg X$, or $X = O(Y)$ to denote the estimate $|X| \leq CY$ for some absolute constant $C$.  
If $x$ is a parameter going to infinity, we use $X = o(Y)$ to denote the claim that $|X| \leq c(x) Y$ for some quantity $c(x)$ that goes to zero as $x \to \infty$ (holding all other parameters fixed).

Unless otherwise specified, all sums are over the integers, except for sums over the variable $p$ (or $p_1$, $p_2$, etc.) which are understood to be over primes.

We use ${\mathbb T} := \R/\Z$ to denote the standard unit circle, and let $e: {\mathbb T} \to \C$ be the standard character $e(x) := e^{2\pi i x}$.

We use $1_S$ to denote the indicator of a predicate $S$, thus $1_S=1$ when $S$ is true and $1_S=0$ when $S$ is false.  If $A$ is a set, we write $1_A(n)$ for $1_{n \in A}$, so that $1_A$ is the indicator function of $A$.

\section{Restricting to numbers with typical factorization}

To prove Theorem \ref{elliott-thm} and Theorem \ref{exp-est} (and hence Theorem \ref{main} and Theorem \ref{main-2}), it is technically convenient (as in the previous paper \cite{MR} of the first two authors) to restrict the support of the multiplicative functions to a certain dense set ${\mathcal S}$ of natural numbers that have a ``typical'' prime factorization in a certain specific sense, in order to fully exploit a useful combinatorial identity of Ramar\'e (see \eqref{ramare} below).  This will lead to improved quantitative estimates in the arguments in subsequent sections of the paper.

More precisely, we introduce the following sets ${\mathcal S}$ of numbers with typical prime factorization, which previously appeared in \cite{MR}.

\begin{definition}\label{S-def}  Let $10 < P_1 < Q_1 \leq X$ and $\sqrt{X} \leq X_0 \leq X$ be quantities such that $Q_1 \leq \exp(\sqrt{\log X_0})$.  We then define
$P_j,Q_j$ for $j > 1$ by the formula
$$ P_j := \exp(j^{4j} (\log Q_1)^{j-1} \log P_1); \quad Q_j := \exp( j^{4j+2} (\log Q_1)^j ).$$
for $j > 1$; note that the intervals $[P_j,Q_j]$ are disjoint and increase to infinity, indeed one easily verifies that 
$$P_1 < Q_1 < \exp( 2^8 \log Q_1 \log P_1 ) = P_2$$ 
and 
$$P_j < \exp( j^{4j} (\log Q_1)^j) < Q_j < \exp( (j+1)^{4(j+1)} (\log Q_1)^j ) < P_{j+1}$$ 
for all $j>1$.  Let $J$ be the largest index such that $Q_J \leq \exp(\sqrt{\log X_0})$.  Then we define ${\mathcal S}_{P_1,Q_1,X_0,X}$ to be the set of all the numbers $1 \leq n \leq X$ which have at least one prime factor in the interval $[P_j,Q_j]$ for each $1 \leq j \leq J$. 
\end{definition}

In practice $X$ will be taken to be slightly smaller than $X_0^2$.  The need to have two parameters $X,X_0$ instead of one is technical (we need to have the freedom later in the argument to replace $X$ with a slightly smaller quantity $X/d$ without altering $J$), but the reader may wish to pretend that $X_0 = \sqrt{X}$ for most of the argument.

This set is fairly dense if $P_1$ and $Q_1$ are widely separated:

\begin{lemma}\label{excep}  Let $10 < P_1 < Q_1 \leq X$ and $\sqrt{X} \leq X_0 \leq X$ be such that $Q_1 \leq \exp( \sqrt{\log X_0} )$.  Then, for every large enough $X$,
$$ \# \{ 1 \leq n \leq X: n \not \in {\mathcal S}_{P_1,Q_1,X_0,X} \} \ll \frac{\log P_1}{\log Q_1} \cdot X.$$
\end{lemma}

\begin{proof}
From the fundamental lemma of sieve theory (see e.g. \cite[Theorem 6.17]{Opera}) we know that, for any $1 \leq j \leq J$ and large enough $X$, the number of $1 \leq n \leq X$ that are not divisible by any prime in $[P_j,Q_j]$ is at most
$$\ll X \prod_{P_j \leq p \leq Q_j} \left(1-\frac{1}{p}\right) \ll \frac{\log P_j}{\log Q_j} X = \frac{1}{j^2} \frac{\log P_1}{\log Q_1} X.$$
Summing over $j$, we obtain the claim.
\end{proof}

Both Theorem \ref{elliott-thm} and Theorem \ref{exp-est} will be deduced from the following claim.

\begin{theorem}[Key exponential sum estimate]\label{second}  Let $X, H, W \geq 10$ be such that
$$ (\log H)^5 \leq W \leq \min\{H^{1/250}, (\log X)^{1/125}\} $$
and let $g$ be a $1$-bounded multiplicative function such that
\begin{equation}\label{wm}
W \leq \exp( M( g; X, Q ) / 3 ).
\end{equation}
Set
$${\mathcal S} := {\mathcal S}_{P_1,Q_1,\sqrt{X},X}$$
where
$$ P_1 := W^{200}; \quad Q_1 := H / W^3.$$
Then for any $\alpha \in {\mathbb T}$, one has
\begin{equation}\label{sumx}
 \int_\R \left|\sum_{x \leq n \leq x+H} 1_{{\mathcal S}}(n) g(n) e(\alpha n)\right|\ dx \ll \frac{(\log H)^{1/4} \log \log H}{W^{1/4}} HX.
\end{equation}
\end{theorem}

In Section \ref{k-tuple-sec} we will show how this theorem implies Theorem \ref{elliott-thm}.  For now, let us at least see how it implies Theorem \ref{exp-est}:

\begin{proof} (Proof of Theorem \ref{exp-est} assuming Theorem \ref{second})  We may assume that $X, H$, and $M(g;X,Q)$ are larger than any specified absolute constant, as if one of these expressions are bounded, then so is $W$, and the claim \eqref{sumx} is then trivial with a suitable choice of implied constant (discarding the $(\log H)^{1/4} \log\log H$ factor).

Choose $H_0$ such that
$$\log H_0 := \min\left( \log^{1/700} X \log\log X, \exp( M(g; X, Q) / 20 ) M(g; X, Q )\right).$$
We divide into two cases: $H \leq H_0$ and $H > H_0$.

First suppose that $H \leq H_0$.  Then if we set $W := \log^5 H$, one verifies that all the hypotheses of Theorem \ref{second} hold, and hence
$$
 \int_0^{X} \left|\sum_{x \leq n \leq x+H} 1_{{\mathcal S}}(n) g(n) e(\alpha n)\right|\ dx \ll \frac{\log \log H}{\log H} HX.
$$
On the other hand, from Lemma \ref{excep}, the choice of $W,P_1,Q_1$, and the bound on $H$ we see that
$$ \# \{ 1 \leq n \leq X+H: n \not \in {\mathcal S} \} \ll \frac{\log\log H}{\log H} X$$
and thus by Fubini's theorem and the triangle inequality
$$ \int_0^{X} \left|\sum_{x \leq n \leq x+H} (1-1_{{\mathcal S}}(n)) g(n) e(\alpha n)\right|\ dx \ll \frac{\log \log H}{\log H} HX.$$
Summing, we obtain Theorem \ref{exp-est} in this case.

Now suppose that $H > H_0$.  Covering $[0,H]$ by $O(H/H_0)$ intervals of length $H_0$, we see that
$$
\int_0^{X} \left|\sum_{x \leq n \leq x+H} g(n) e(\alpha n)\right|\ dx \ll 
\frac{H}{H_0}
\int_{0}^{X+H} \left|\sum_{x \leq n \leq x+H_0} g(n) e(\alpha n)\right|\ dx.
$$
Also, observe from the choice of $H_0$ that the quantity $\exp( -M(g; X, Q) / 20 ) + \frac{\log \log H}{\log H} + \frac{1}{\log^{1/700} X}$ is unchanged up to multiplicative constants if one reduces $H$ to $H_0$.  Finally, from Mertens' theorem we see that $M(g;X+H,Q) = M(g;X,Q) + O(1)$.  The claim then follows from the $H=H_0$ case (after performing the minor alteration of replacing $X$ with $X+H$).
\end{proof}

We now begin the proof of Theorem \ref{second}.  The first step is to reduce to the case where $g$ is completely multiplicative rather than multiplicative.  More precisely, we will deduce Theorem \ref{second} from

\begin{proposition}[Completely multiplicative exponential sum estimate]\label{third}  Let $X, H, W \geq 10$ be such that
$$ (\log H)^{5} \leq W \leq \min\{H^{1/250}, (\log X)^{1/125}\},$$
and let $g$ be a $1$-bounded completely multiplicative function such that
\begin{equation}\label{woo}
 W \leq \exp( M( g; X, W ) / 3 ).
\end{equation}
Let $d$ be a natural number with $d<W$.  Set
$${\mathcal S} := {\mathcal S}_{P_1,Q_1,\sqrt{X},X/d}$$
where
$$ P_1 := W^{200}; \quad Q_1 := H / W^3.$$
Then for any $\alpha \in {\mathbb T}$ one has
\begin{equation}\label{sumx-new}
 \int_\R \left|\sum_{x/d \leq n \leq x/d+H/d} 1_{{\mathcal S}}(n) g(n) e(\alpha n)\right|\ dx \ll \frac{1}{d^{3/4}} \frac{(\log H)^{1/4} \log \log H}{W^{1/4}} HX.
\end{equation}
\end{proposition}

Let us explain why Theorem \ref{second} follows from Proposition \ref{third}.  Let the hypotheses and notation be as in Theorem \ref{second}.  The function $g$ is not necessarily completely multiplicative, but we may approximate it by the $1$-bounded completely multiplicative function $g_1: \N \to \C$, defined as the completely multiplicative function with $g_1(p)=g(p)$ for all primes $p$.  By M\"obius inversion we may then write $g = g_1 * h$ where $*$ denotes Dirichlet convolution and $h$ is the multiplicative function $h = g * \mu g_1$.  Observe that for all primes $p$, $h(p)=0$ and $|h(p^j)| \leq 2$ for $j \geq 2$.  We now write
$$ \sum_{x \leq n \leq x+H} 1_{{\mathcal S}_{P_1,Q_1, \sqrt{X}, X}} g(n) e(\alpha n) = \sum_{d=1}^\infty h(d) \sum_{x/d \leq m \leq x/d + H/d} 1_{{\mathcal S}_{P_1,Q_1, \sqrt{X}, X}}(dm) g_1(m) e(d\alpha m)$$
and so by the triangle inequality we may upper bound the left-hand side of \eqref{sumx} by
$$
\sum_{d=1}^\infty |h(d)| \int_{\mathbb{R}} \left|\sum_{x/d \leq m \leq x/d + H/d} 1_{{\mathcal S}_{P_1,Q_1, \sqrt{X}, X}}(dm) g_1(m) e(d\alpha m)\right| dx.$$
Let us first dispose of the contribution where $d \geq W$.  Here we trivially bound this contribution by
$$ \sum_{d \geq W} |h(d)| \sum_{m \leq (2X+H)/d} O(H) $$
(after moving the absolute values inside the $m$ summation and then performing the integration on $x$ first).  We can bound this in turn by
$$ \ll H X \frac{1}{W^{1/4}} \sum_{d=1}^\infty \frac{|h(d)|}{d^{3/4}}.$$
From Euler products we see that $\sum_{d=1}^\infty \frac{|h(d)|}{d^{3/4}} = O(1)$, so the contribution of this case is acceptable.

Now we consider the contribution $d < W < P_1$.  In this case we may reduce
$$ 1_{{\mathcal S}_{P_1,Q_1,\sqrt{X},X}}(dm) = 1_{{\mathcal S}_{P_1,Q_1,\sqrt{X},X/d}}(m)$$
and so this contribution to \eqref{sumx} can be upper bounded by
$$
\sum_{1 \leq d < W} |h(d)| \int_\R \left|\sum_{x/d \leq m \leq x/d + H/d} 1_{{\mathcal S}_{P_1,Q_1,\sqrt{X},X/d}}(m) g_1(m) e(d\alpha m)\right|\ dx.$$
By Proposition \ref{third}, this is bounded by
$$ \sum_{d=1}^\infty \frac{|h(d)|}{d^{3/4}} \frac{(\log H)^{1/4} \log \log H}{W^{1/4}} HX.$$
As before we have $\sum_{d=1}^\infty \frac{|h(d)|}{d^{3/4}} = O(1)$, and Theorem \ref{second} follows.

It remains to prove Proposition \ref{third}.  For any $\alpha \in {\mathbb T}$, we know from the Dirichlet approximation theorem that there exists a rational number $\frac{a}{q}$ with $(a,q)=1$ and $1 \leq q \leq H/W$ such that
$$ \left|\alpha - \frac{a}{q}\right| \leq \frac{W}{qH} \leq \frac{1}{q^2}.$$
In the next two sections, we will apply separate arguments to prove Proposition \ref{third} in the \emph{minor arc} case $q > W$ and the \emph{major arc} case $q \leq W$.

\section{Proof of minor arc estimate}

We now prove Proposition \ref{third} in the minor arc case $q > W$.  It suffices to show that
\begin{equation}\label{mainsum} 
 \int_\R \theta(x) \sum_{x/d \leq n \leq x/d+H/d} 1_{{\mathcal S}}(n) g(n) e(\alpha n)\ dx \ll
\frac{1}{d^{3/4}} \frac{(\log H)^{1/4} \log \log H}{W^{1/4}} HX
\end{equation}
whenever $\theta: \R \to \C$ is measurable with $|\theta(x)| \leq 1$ for all $x$ and supported on $[0, X]$. 
We will now use a variant of an idea of Bourgain-Sarnak-Ziegler \cite{BSZ} (building on earlier works of 
Katai \cite{Katai}, Montgomery-Vaughan \cite{MontgomeryVaughan} and Daboussi-Delange \cite{DaboussiDelange}).

Let $\mathcal{P}$ be the set consisting of the primes lying between $P_1$ and $Q_1$. Then, notice that each $n \in \mathcal{S}$ has at least one prime factor from $\mathcal{P}$.  Furthermore, if $n=mp$ for some prime $p \in \mathcal{P}$, then the number of primes in $\mathcal{P}$ dividing $n$ is equal to the number of primes in $\mathcal{P}$ dividing $m$, plus\footnote{In the published version of this paper, the term $1_{p \nmid m}$ was incorrectly expressed as $1$, leading to a slight gap in the arguments.  We thank Alisa Sedunova and Ke Wang for drawing this issue to our attention.} $1_{p \nmid m}$.  This leads to the following variant of Ramar\'e's identity (see \cite[Section 17.3]{Opera}):
\begin{equation}\label{ramare}
1_{{\mathcal S}}(n) =  \sum_{p \in {\mathcal P}, m: mp = n} \frac{1_{{\mathcal S}'}(mp)}{1_{p \nmid m} + \# \{ q \mid m: q \in {\mathcal P}\}},
\end{equation}
where ${\mathcal S}'$ is the set of all $1 \leq n \leq X/d$ that have at least one prime factor in each of the intervals $[P_j,Q_j]$ for $j \geq 2$; the constraint $n \leq X/d$ arises from the corresponding constraint in the definition of ${\mathcal S}$.

Using this identity, we may write the left-hand side of \eqref{mainsum} as
$$
\sum_{p \in \mathcal{P}} \sum_{m} \frac{1_{{\mathcal S}'}(m p) g(m p) e (m p \alpha)}{1_{p \nmid m} + \#\{q | m: q \in \mathcal{P}\}} \int_\R \theta(x) 1_{x/d \leq mp \leq (x+H)/d}\ dx.
$$
As $g$ is completely multiplicative, $g(mp) = g(m) g(p)$.  Thus it suffices to show that
\begin{align*}
\sum_{p \in \mathcal{P}} \sum_{m}  \frac{1_{{\mathcal S}'}(m p) g(m) g(p) e (m p \alpha)}{1_{p \nmid m} + \#\{q | m: q \in \mathcal{P}\}} \int_{\mathbb{R}} & \theta(x) 1_{x/d \leq mp \leq (x+H)/d} \ dx \ll \\ &
\ll 
\frac{(\log H)^{1/4} \log \log H}{d^{3/4} W^{1/4}} HX.
\end{align*}
We can cover ${\mathcal P}$ by intervals $[P,2P]$ with $P_1 \ll P \ll Q_1$ and $P$ a power of two, and observe that
$$ \sum_{\substack{P_1 \ll P \ll Q_1 \\ P = 2^j}} \frac{1}{\log P} \ll \log \log Q_1 - \log\log P_1 \ll \log\log H,$$
so by the triangle inequality it suffices to show that
$$
\sum_{p \in {\mathcal P}: P \leq p \leq 2P} \sum_{m} \frac{1_{{\mathcal S}'}(m p) g(m) g(p) e (m p \alpha)}{1_{p \nmid m} + \#\{q | m: q \in \mathcal{P}\}} \int_{\mathbb{R}} \theta(x) 1_{x/d \leq mp \leq (x+H)/d} \ dx
\ll 
\frac{(\log H)^{1/4}}{d^{3/4} W^{1/4} \log P} HX
$$
for each such $P$.  

Fix $P$.  At this point it becomes convenient to replace the $1_{p \nmid m}$ term by $1$.  Since the integral $\int_{\mathbb{R}} \theta(x) 1_{x/d \leq mp \leq (x+H)/d} \ dx$ is $O(H)$ by the triangle inequality, and all the other factors in the summand are $O(1)$, and the term $1_{{\mathcal S}'}(mp)$ vanishes unless $m \leq X/dP$, the error incurred in making this substitution may be bounded in magnitude by
$$ O \left( \sum_{p \in {\mathcal P}: P \leq p \leq 2P}  \sum_{m \leq X/dP: p|m} H \right) = O\left( P \frac{X}{dP^2} H \right) = O\left( \frac{HX}{dP} \right),$$
 which is acceptable.  Thus it will suffice to show that
$$
\sum_{p \in {\mathcal P}: P \leq p \leq 2P} \sum_{m} \frac{1_{{\mathcal S}'}(m p) g(m) g(p) e (m p \alpha)}{1 + \#\{q | m: q \in \mathcal{P}\}} \int_{\mathbb{R}} \theta(x) 1_{x/d \leq mp \leq (x+H)/d} \ dx
\ll 
\frac{(\log H)^{1/4}}{d^{3/4} W^{1/4} \log P} HX.
$$
We can rearrange the left-hand side as
$$
\sum_{m \in {\mathcal S}'} \frac{g(m)}{1 + \#\{q | m: q \in \mathcal{P}\}}
\sum_{p \in {\mathcal P}: P \leq p \leq 2P} 1_{mp \leq X/d} g(p) e (m p \alpha) \int_\R \theta(x) 1_{x/d \leq mp \leq (x+H)/d}\ dx.$$
As before, the summand vanishes unless $m \leq \frac{X}{dP}$.  Crudely bounding\footnote{By using the Turan-Kubilius inequality here one could save a factor of $\log\log H$, but such a gain will not make a significant impact on our final estimates.} $\frac{g(m)}{1 + \#\{q | m: q \in \mathcal{P}\}}$ in magnitude by $1$, we may bound the previous expression in magnitude by
$$
\sum_{m \leq X/dP} \left|
\sum_{p \in {\mathcal P}: P \leq p \leq 2P} 1_{mp \leq X/d} g(p) e (m p \alpha) \int_\R \theta(x) 1_{x/d \leq mp \leq (x+H)/d}\ dx\right|.$$
By H\"older, we may bound this by
$$ \left( \frac{X}{dP} \right)^{3/4}
\left( 
\sum_{m \leq X/dP}
\left|\sum_{p \in {\mathcal P}: P \leq p \leq 2P} 1_{mp \leq X/d} g(p) e (m p \alpha) \int_{\mathbb{R}} \theta(x) 1_{x/d \leq mp \leq (x+H)/d}\right|^4 dx \right)^{1/4}.$$
It thus suffices to show that
$$
\sum_{m \leq X/dP}
\left|\sum_{p \in {\mathcal P}: P \leq p \leq 2P} 1_{mp \leq X/d} g(p) e (m p \alpha) \int_{\mathbb{R}} \theta(x) 1_{x/d \leq mp \leq (x+H)/d} dx\right|^4 
\ll
\frac{\log H}{W \log^4 P} H^4 X P^3.$$
The left-hand side may be expanded as
\begin{align*}
&\sum_{p_1,p_2,p_3,p_4 \in {\mathcal P}: P \leq p_1,p_2,p_3,p_4 \leq 2P} 
\int \dotsi \int g(p_1) g(p_2) \overline{g(p_3)} \overline{g(p_4)} \theta(x_1) \theta(x_2) \overline{\theta(x_3)} \overline{\theta(x_4)} \\
&\cdot \sum_{m \leq X/(dp_i), x_i/(dp_i) \leq m \leq (x_i+H)/(dp_i) \forall i=1,2,3,4} e(m(p_1+p_2-p_3-p_4) \alpha) dx_1 dx_2 dx_3 dx_4.
\end{align*}
From summing the geometric series, the summation over $m$ is $O( \min( \frac{H}{P}, \frac{1}{\|(p_1+p_2-p_3-p_4) \alpha\|} ) )$, where $\|z\|$ denotes the distance from $z$ to the nearest integer.  Also, the sum vanishes unless we have $x_1 = O(X)$ and $x_i = x_1p_i/p_1 + O( H)$ for $i=2,3,4$, so there are only $O( X H^3 )$ quadruples $(x_1,x_2,x_3,x_4)$ which contribute here.  Thus we may bound the previous expression by
$$
O\left( X H^3 \sum_{p_1,p_2,p_3,p_4 \leq 2P} \min\left(\frac{H}{P}, \frac{1}{\|(p_1+p_2-p_3-p_4) \alpha\|}\right) \right)
$$
and so we reduce to showing that
\begin{equation}\label{flask}
 \sum_{p_1,p_2,p_3,p_4 \leq 2P} \min\left(\frac{H}{P}, \frac{1}{\|(p_1+p_2-p_3-p_4) \alpha\|} \right) \ll \log H \frac{H P^3}{W \log^4 P}.
\end{equation}
The quantity $p_1+p_2-p_3-p_4$ is clearly of size $O(P)$. Conversely, from a standard upper bound sieve\footnote{For instance, from \cite[Theorem 3.13]{MontgomeryVaughanBook} one sees that any number $N = O(P)$ has $O( \frac{N}{\phi(N)} \frac{P}{\log^2 P} )$ representations as the sum of two primes; since $\sum_{N=O(P)} \frac{N^2}{\phi(N)^2} = O(P)$ (see e.g. \cite[Exercise 2.1.14]{MontgomeryVaughanBook}), the claim then follows from the Cauchy-Schwarz inequality.}, the number of representations of an integer $n = O(P)$ of the form $p_1+p_2-p_3-p_4$ with $p_1,p_2,p_3,p_4 \leq 2P$ prime is $O( \frac{P^3}{\log^4 P} )$.  Thus it suffices to show that
$$
 \sum_{n = O(P)} \min\left(\frac{H}{P}, \frac{1}{\|n \alpha\|} \right) \ll \frac{\log H}{W} H.$$
But from the Vinogradov lemma (see e.g. \cite[Page 346]{IwKo04}), the left-hand side is bounded by
$$ O\left( \left(\frac{P}{q}+1\right) \left(\frac{H}{P} + q \log q\right) \right) \ll \frac{H}{q} + P \log q + \frac{H}{P} + q\log q$$
which, since 
$$ W^{200} = P_1 \ll P \ll Q_1 = H/W^3$$
and
$$ W \leq q \leq H/W,$$
is bounded by $O( \frac{\log H}{W} H )$ as required.

\section{Proof of major arc estimate}

We now prove Proposition \ref{third} in the major arc case $q \leq W$.
We will discard the factor $d^{1/4} (\log H)^{1/4} \log\log H$ and prove the 
following stronger bound
\begin{equation} \label{equ:statement}
\int_{\mathbb{R}} \Big | \sum_{x / d \leq n \leq (x + H) / d}
1_{\mathcal{S}}(n) g(n) e(\alpha n) \Big | dx \ll
\frac{H X}{d W^{1/4}}. 
\end{equation}
By hypothesis we have $\alpha = a/q + \theta$ 
with $q \leq W$ and
$\theta = O(W / (H q))$. Integrating by parts we see
that
\begin{align} \nonumber
\Big | \sum_{x /d \leq n \leq (x + H)/ d} 1_{\mathcal{S}}(n) g(n) e(\alpha n)
\Big | & \ll \Big | \sum_{x/d \leq n \leq (x + H)/ d} 1_{\mathcal{S}}(n)
g(n) e(a n / q) \Big | \\ \label{equ:parts}
& + \frac{W}{H q} \int_{0}^{H/d} \Big | \sum_{x / d \leq n \leq x / d + H'}
1_{\mathcal{S}}(n) g(n) e(a n / q) \Big | d H' 
\end{align}
Thus let us focus on bounding,
\begin{equation} \label{equ:goal}
\int_{\mathbb{R}} \Big | \sum_{x/d \leq n \leq x/d + H'} 1_{\mathcal{S}}(n) g(n)
e(a n / q) \Big | dx 
\end{equation}
with $0 \leq H' \leq H / d$. Splitting into residues classes we see that
(\ref{equ:goal}) is
$$
\leq \sum_{b \pmod{q} } \int_{\mathbb{R}} \Big | \sum_{\substack{x /d \leq n \leq x / d + H' \\ n \equiv b \pmod{q}}} 1_{\mathcal{S}}(n) g(n) \Big | dx
$$
For $n \equiv b \pmod {q}$ we have $d_0 := (b,q) | n$. Therefore
let us write $b = d_0 b_0, q = d_0 q_0$ and $n = d_0 m $, so that
the condition $n \equiv b \pmod{q}$ simplifies to $m \equiv b_0 \pmod{q_0}$. 
In addition, 
since $g$ is completely multiplicative and $d_0 \leq q \leq W \leq P_1$
we have
$$
1_{\mathcal{S}}(n) g(n) = g(d_0) \cdot 1_{\mathcal{S}_{P_1, Q_1, \sqrt{X}, X / (d d_0)}} (m) g(m).
$$
Finally we express $m \equiv b_0 \pmod{q_0}$ in terms of Dirichlet characters
noting that
$$
1_{m \equiv b_0 \pmod{q_0}}(m) = \frac{1}{\varphi(q_0)} 
\sum_{\chi \pmod{q_0}} \chi(b_0) \overline{\chi(m)}. 
$$
Plugging everything together we see that (\ref{equ:goal}) is less than
$$
\sum_{b \pmod{q}} \frac{1}{\varphi(q_0)} \sum_{\chi \pmod{q_0}}
\int_{\mathbb{R}} \Big | \sum_{x / (d d_0) \leq m \leq x / (d d_0) + H'/d_0}
1_{\mathcal{S}_{P_1, Q_1, \sqrt{X}, X / (d d_0)}} g(m) \overline{\chi(m)}
\Big | dx
$$
In the integral we make the linear change of variable $y = x / (d d_0)$, so that
the above expression becomes
\begin{equation} \label{equ:goal2}
d \sum_{b \pmod{q}} \frac{d_0}{\varphi(q_0)} \sum_{\chi \pmod{q_0}} 
\int_{\mathbb{R}} \Big | \sum_{y \leq m \leq y + H'/d_0} 1_{\mathcal{S}_{P_1, Q_1, \sqrt{X}, X / (d d_0)}} g(m) \overline{\chi(m)} \Big | dy 
\end{equation}
We bound the part of the integral with $y \leq X / W^{10}$ trivially. This produces in (\ref{equ:goal}) an error which is
$$
\ll d q \cdot \frac{X}{W^{10}} \cdot H'  \leq \frac{H X}{W^{9}} \ll \frac{H X}{d W^{3}}
$$
since $q , d \leq W$ and $H' \leq H / d$. We split the remaining range $X / W^{10} \leq y \leq 2X / (d d_0)$ into dyadic blocks $X / W^{10} \leq X' \leq X / (d d_0)$ with $X'$ running through powers of two. Thus the previous expression is
\begin{align*}
\ll d \sum_{X'} & \sum_{b \pmod{q}} \frac{d_0}{\varphi(q_0)} \sum_{\chi \pmod{q_0}} 
\int_{X'}^{2X'} \Big | \sum_{y \leq m \leq y + H' /d_0} 1_{\mathcal{S}_{P_1, Q_1, \sqrt{X}, X / (d d_0)}} g(m)
\overline{\chi(m)} \Big | dy \\ & + \frac{H X}{d W^3}. 
\end{align*}
At this point we apply Theorem \ref{th:MRinC} with $\eta=1/20$ (note that $P_1 \geq (\log Q_1)^{40/\eta}$) to conclude that
\[                                                                                                                                   
\begin{split}                                                                                                                          
&\int_{X'}^{2X'} \left|\sum_{y \leq m \leq y+H'/d_0} 1_{{\mathcal S}_{P_1,Q_1,\sqrt{X},X/(d d_0)}}(m) g(m) \overline{\chi}(m)\right|^2\ dy \\ 
&\ll \left(\exp( -M( g \overline{\chi}; X') ) M(g \overline{\chi}; X') + \frac{(\log H'/d_0)^{1/3}}{P_1^{1/6-1/20}} + \frac{1}{(\log X')^{1/\
50}} \right) \frac{H'^2}{d_0^2} X'.                                                                                                                 
\end{split}                                                                                                                             
\]
Since $P_1 = W^{200}$ and $H'/d_0 \leq H$ and $W \geq \log^{5} H$, we have
$$ \frac{(\log H'/d_0)^{1/3}}{P_1^{1/6 - 1/20}} \leq \frac{(\log H)^{1/3}}{P_1^{1/6-1/20}} \ll \frac{1}{W^{5/2}} $$
and certainly
$$ \frac{1}{(\log X')^{1/50}} \ll \frac{1}{(\log X)^{1/50}} \ll \frac{1}{W^{5/2}}.$$
From Mertens' theorem and definition of $M(g, X, W)$,
$$ M( g \overline{\chi}; X') \geq M( g \overline{\chi}; X) - O(1) \geq M(g, X, W) - O(1)$$
and thus by \eqref{woo}
$$ \exp( -M( g \overline{\chi}; X') ) M(g \overline{\chi}; X') \ll \frac{1}{W^{5/2}}.$$
Putting all this together, we obtain
$$ \int_{X'}^{2X'} \left|\sum_{y \leq m \leq y+H'/d_0} 1_{{\mathcal S}_{P_1,Q_1,\sqrt{X},X/(d d_0)}}(m) g(m) \overline{\chi}(m)\right|^2\ dy
\ll \frac{1}{W^{5/2}} \frac{H'^2}{d_0^2} X'.$$
It follows from Cauchy-Schwarz that 
$$
\int_{X'}^{2X'} \Big | \sum_{y \leq m \leq y + H'/d_0} 1_{\mathcal{S}_{P_1, Q_1, \sqrt{X}, X/(d d_0)}}(m) g(m)
\overline{\chi(m)} \Big | dy \ll W^{-5/4} \cdot \frac{H' X'}{d_0} 
$$
Inserting this bound into (\ref{equ:goal2}) we see that (\ref{equ:goal}) is bounded by
$$
\ll d q \cdot \frac{1}{W^{5/4}} \cdot \frac{H}{d} \cdot \frac{X}{d} \ll \frac{q H X}{d W^{5/4}} 
$$
Therefore using (\ref{equ:parts}) and using $q \leq W$ we see that (\ref{equ:statement}) is
$$
\ll \frac{q H X}{d W^{5/4}} \cdot \Big (1 + \frac{W}{H q} \cdot \frac{H}{d} \Big )
\ll \frac{H X}{d W^{1/4}}
$$
as claimed. 

\section{Elliott's conjecture on the average}\label{k-tuple-sec}

In this section we use Theorem \ref{second} to prove Theorem \ref{elliott-thm}. Theorem~\ref{elliott-thm} will be deduced from the following result (compare also with Theorem~\ref{second} and deduction of Theorem~\ref{exp-est} from it).  For brevity, we write $1_{\mathcal S} g$ for the function $n \mapsto 1_{\mathcal S}(n) g(n)$.

\begin{proposition}[Truncated Elliott on the average]\label{first}  
Let $X, H, W, A \geq 10$ be such that
\[
\log^{20} H \leq W \leq  \min\{H^{1/500}, (\log X)^{1/125}\}.
\]

Let $g_1,\dots,g_k: \N \to \C$ be $1$-bounded multiplicative functions, and let $a_1,\dots,a_k,b_1,\dots,b_k$ be natural numbers with $a_j \leq A$ and $b_j \leq 3AX$ for $j=1,\dots,k$.  Let $1 \leq j_0 \leq k$ be such that
$$ W \leq \exp(M( g_{j_0}; 10AX, Q ) / 3 ).$$

Set
$$ {\mathcal S} = {\mathcal S}_{P_1, Q_1, \sqrt{10AX}, 10AX}$$
where
$$ P_1 := W^{200}; \quad Q_1 := H^{1/2} / W^3.$$ 
Then
\begin{equation}\label{first-eq}
 \sum_{1 \leq h_2,\dots,h_k \leq H} \left|\sum_{1 \leq n \leq X} 1_{\mathcal S} g_1(a_1 n + b_1) \prod_{j=2}^k 1_{\mathcal S} g_j(a_j n + b_j + h_j)\right| \ll \frac{kA^2}{W^{1/20}} H^{k-1} X.
\end{equation}
\end{proposition}

\begin{proof}[Proof of Theorem~\ref{elliott-thm} assuming Proposition~\ref{first}]
We may assume that $X$, $H$, and $M$ are larger than any specified absolute constant as the claim is trivial otherwise. We first make some initial reductions.  The first estimate \eqref{ag} of Theorem \ref{elliott-thm} follows from the second \eqref{bg} after shifting $b_1$ by $h_1$ in \eqref{bg} and averaging, provided that we relax the hypotheses $b_j \leq AX$ slightly to $b_j \leq 2AX$.  Thus it suffices to prove \eqref{bg} under the relaxed hypotheses $b_j \leq 2AX$.

Let $H_0$ be such that 
\begin{equation}\label{hax}
 \log H_0 = \min\{\log^{1/3000} X \log\log X, \exp(M(g_{j_0};10AX, Q)/80) M(g_{j_0};10AX, Q)\}. 
\end{equation}
If $H \leq H_0$ we take $W = \log^{20} H$ and let $\mathcal{S}$ be as in Proposition~\ref{first}. All the assumptions of Proposition~\ref{first} hold and thus
\[
 \sum_{1 \leq h_2,\dots,h_k \leq H} \left|\sum_{1 \leq n \leq X} 1_{\mathcal S} g_1(a_1 n + b_1) \prod_{j=2}^k 1_{\mathcal S} g_j(a_j n + b_j + h_j)\right| \ll \frac{kA^2}{\log H} H^{k-1} X.
\]

Furthermore, from Lemma \ref{excep} we have
\begin{equation}\label{sdense}
\sum_{n \leq 10AX: n \not \in {\mathcal S}} 1 \ll AX \frac{\log W}{\log H}.
\end{equation}
From this and the triangle inequality, we have
\begin{equation}\label{ggg}
\begin{split}
\sum_{1 \leq n \leq X} & g_1(a_1 n + b_1)  \prod_{j=2}^k g_j(a_j n + b_j + h_j) \\
& = \sum_{1 \leq n \leq X} 1_{\mathcal S} g_1(a_1 n + b_1) \prod_{j=2}^k 1_{\mathcal S} g_j(a_j n + b_j + h_j) + O \Big ( k A X \frac{\log W}{\log H} \Big ).
\end{split}
\end{equation}
Hence the claim follows in the case when $H \leq H_0$.

If $H > H_0$, one can cover the summation over the $h_j$ indices by intervals of length $H_0$ and apply Theorem \ref{elliott-thm} to each subinterval (shifting the $b_j$ by at most $AX$ when doing so), and then sum, noting that the quantity 
$$ \exp(-M(g_{j_0};10AX, Q )/80) + \frac{\log\log H}{\log H} + \frac{1}{\log^{1/3000} X}$$ 
is essentially unchanged after replacing $H$ with $H_0$.
\end{proof}

\begin{remark}\label{improv}  By using larger choices of $W$, one can obtain more refined information on the large values of the correlations $\sum_{1 \leq n \leq X} g_1(a_1 n + b_1) \prod_{j=2}^k g_j(a_j n + b_j + h_j)$.  For instance, if we take $W = H^{\delta}$ for some $10 \leq H \leq H_0$ and $20 \frac{\log\log H}{\log H} \leq \delta \leq \frac{1}{500}$, we see from Proposition \ref{first}, \eqref{ggg}, and Markov's inequality that
$$
\sum_{1 \leq n \leq X} g_1(a_1 n + b_1) \prod_{j=2}^k g_j(a_j n + b_j + h_j) \ll k A^2 \delta X
$$
for all but at most $O( \frac{H^{k-1}}{\delta H^{\delta/20}} )$ tuples $(h_1,\dots,h_{k-1})$ with $1 \leq h_j \leq H$ for $j=2,\dots,k$.  Thus we can obtain a power saving in the number of exceptional tuples, at the cost of only obtaining a weak bound on the individual correlations
$\sum_{1 \leq n \leq X} g_1(a_1 n + b_1) \prod_{j=2}^k g_j(a_j n + b_j + h_j)$.  
\end{remark}

It remains to prove Proposition \ref{first}. We start by proving the following simpler case to which the general case will be reduced.
\begin{proposition}\label{prop:truncEllsc}
Let $X, H, W \geq 10$ be such that
\[
\log^{20} H \leq W \leq  \min\{H^{1/250}, (\log X)^{1/125}\}.
\]
Let $g: \N \to \C$ be $1$-bounded multiplicative function such that
$$ W \leq \exp(M( g; X, W ) / 3 ).$$

Set
$$ {\mathcal S} = {\mathcal S}_{P_1, Q_1, \sqrt{X}, X}$$
where
$$ P_1 := W^{200}; \quad Q_1 := H/ W^3.$$ 
Then
\begin{equation}\label{first-eq2}
 \sum_{1 \leq h \leq H} \left|\sum_{1 \leq n \leq X} 1_{\mathcal S} g(n) 1_{\mathcal S}\overline{g}(n+h)\right|^2 \ll \frac{H X^2}{W^{1/5}}.
\end{equation}
\end{proposition}

To deduce Theorem~\ref{th:quant} we let $\mathcal{S}$ be as in this proposition with $W := H^{\delta/900}$. The argument of Lemma~\ref{excep} actually gives in this case $\#\{1 \leq n \leq X \colon n \not \in \mathcal{S}\} \ll \frac{\log P_1}{\log Q_1} X + \frac{X}{P_1}$, and thus the numbers $n$ with $n \not \in \mathcal{S}$ or $n+h \not \in \mathcal{S}$ contribute to the left hand side of~\eqref{eq:quanteq} at most $9\delta/10$. Hence, recalling~\eqref{eq:MestLiou}, the claim follows from the previous proposition and Markov's inequality.

\begin{proof}[Proof of Proposition~\ref{prop:truncEllsc}]
The claim follows once we have shown
$$ \sum_{|h| \leq 2H} (2H-|h|)^2 \cdot \Big |\sum_n 1_{\mathcal S} g(n) 1_{\mathcal S} \overline{g}(n+h) \Big |^2 \ll \frac{1}{W^{1/5}} H^3 X^2.$$
Applying Lemma \ref{fourier}, it will suffice to show that
$$ \int_{{\mathbb T}} \left( \int_\R \left|\sum_{x \leq n \leq x+2H} 1_{\mathcal S} g(n) e(\alpha n)\right|^2\ dx \right)^2\ d\alpha \ll \frac{1}{W^{1/5}} H^3 X^2.$$
From the Parseval identity we have
\begin{align*}
\int_{{\mathbb T}} \int_\R \left|\sum_{x \leq n \leq x+2H} 1_{\mathcal S} g(n) e(\alpha n)\right|^2\ dx \ d\alpha\
&= \int_\R \sum_{x \leq n \leq x+2H} |1_{\mathcal S} g(n)|^2\ dx \\
&\ll H  X 
\end{align*}
so it suffices to show that
$$ \sup_\alpha \int_\R \Big |\sum_{x \leq n \leq x+2H} 1_{\mathcal S} g(n) e(\alpha n) \Big |^2\ dx \ll \frac{1}{ W^{1/5}} H^2 X.$$
Using the trivial bound
$$ \Big |\sum_{x \leq n \leq x+2H} 1_{\mathcal S} g(n) e(\alpha n) \Big | \ll H$$
we thus reduce to showing 
\begin{equation}\label{war}
 \sup_\alpha \int_\R \Big |\sum_{x \leq n \leq x+2H} 1_{\mathcal S} g(n) e(\alpha n) \Big |\ dx \ll \frac{H X}{ W^{1/5}}.
\end{equation}
But this follows from Theorem \ref{second} (using the lower bound $W \geq \log^{20} H$ in the hypotheses of Proposition~\ref{prop:truncEllsc} to absorb the $\log^{1/4} H \log\log H$ factors in Theorem \ref{second}).
\end{proof}

\begin{proof}[Proof of Proposition~\ref{first}]
We first remove the special treatment afforded to the $g_1$ factor in \eqref{first-eq}.  Note that we may assume that
\begin{equation}\label{wham}
 W^{1/20} \geq kA^2
\end{equation}
and thus
$$ H \geq W^{500} \geq (kA^2)^{10000}$$
since the claim is trivial otherwise.

Set $H' := \sqrt{H}$.  For any $1 \leq h_1 \leq H'/A$, we may shift $n$ by $h_1$ and conclude that
\begin{align*}
&\sum_{1 \leq n \leq X} 1_{\mathcal S} g_1(a_1 n + b_1) \prod_{j=2}^k 1_{\mathcal S} g_j(a_j n + b_j + h_j) \\
&= 
\sum_{1 \leq n \leq X} 1_{\mathcal S} g_1(a_1 n + b_1 + a_1 h_1) \prod_{j=2}^k 1_{\mathcal S} g_j(a_j n + b_j + h_j + a_j h_1) + O( H' )
\end{align*}
and thus we may write the left-hand side of \eqref{first-eq} as
$$
 \sum_{1 \leq h_2,\dots,h_k \leq H} \left|\sum_{1 \leq n \leq X} 1_{\mathcal S}g_1(a_1 n + b_1 + a_1 h_1) \prod_{j=2}^k 1_{\mathcal S}g_j(a_j n + b_j + h_j + a_j h_1)\right| + O( H^{k-1} H' ).
$$
If one shifts each of the $h_j$ for $j=2,\dots,k$ in turn by $a_j h_1 = O( H' )$, we may rewrite this as
$$
 \sum_{1 \leq h_2,\dots,h_k \leq H} \left|\sum_{1 \leq n \leq X} 1_{\mathcal S}g_1(a_1 n + b_1 + a_1 h_1) \prod_{j=2}^k 1_{\mathcal S}g_j(a_j n + b_j + h_j)\right| + O( H^{k-1} H' ) + O( k H^{k-2} H' X ).
$$
Averaging in $h_1$, and replacing $h_1$ by $a_1 h_1$ (crudely dropping the constraint that $a_1 h_1$ is divisible by $a_1$), we may thus bound the left-hand side of \eqref{first-eq} by
\begin{align*} \ll \frac{A}{H'} \sum_{1 \leq h_1 \leq H'}
 \sum_{1 \leq h_2,\dots,h_k \leq H} & \left|\sum_{1 \leq n \leq X} 1_{\mathcal S}g_1(a_1 n + b_1 + h_1) \prod_{j=2}^k 1_{\mathcal S}g_j(a_j n + b_j + h_j)\right| + \\ & + H^{k-1}H'+ k H^{k-2}H' X.
\end{align*}
The $g_1$ term may now be combined with the product over the remaining $g_j$ terms to form $\prod_{j=1}^k 1_{\mathcal S}g_j(a_j n + b_j + h_j)$.  The error term $H^{k-1}H'+ k H^{k-2}H' X$ is certainly of size $O( \frac{kA^2}{W^{1/20}} H^{k-1} X )$, so it suffices to show that
$$
\sum_{1 \leq h_1 \leq H'}
 \sum_{1 \leq h_2,\dots,h_k \leq H} \left|\sum_{1 \leq n \leq X} \prod_{j=1}^k 1_{\mathcal S}g_j(a_j n + b_j + h_j)\right| \ll \frac{A}{W^{1/20}} H^{k-1}H' X.
$$
By covering the ranges $1 \leq h_j \leq H$ by intervals of length $H'$ and averaging, it suffices (after relaxing the conditions $b_j \leq 3AX$ to $b_j \leq 4AX$) to prove that
$$
\sum_{1 \leq h_1,h_2,\dots,h_k \leq H'} \left|\sum_{1 \leq n \leq X} \prod_{j=1}^k 1_{\mathcal S}g_j(a_j n + b_j + h_j)\right| \ll \frac{A}{W^{1/20}} (H')^k X.
$$
The situation is now symmetric with respect to permuting the indices $1,\dots,k$, so we may assume that the index $j_0$ in Proposition \ref{first} is equal to $1$.  By the triangle inequality in $h_2,\dots,h_k$, it suffices to show that
$$
\sum_{1 \leq h_1 \leq H'} \left|\sum_{1 \leq n \leq X} \prod_{j=1}^k 1_{\mathcal S}g_j(a_j n + b_j + h_j)\right| \ll \frac{A}{W^{1/20}} H' X
$$
for all $h_2,\dots,h_k$.  Writing $G(n) := \prod_{j=2}^k 1_{\mathcal S}g_j(a_j n + b_j + h_j)$, it thus suffices to show that
$$
\sum_{1 \leq h_1 \leq H'} \left|\sum_{1 \leq n \leq X} 1_{\mathcal S}g_1(a_1 n + b_1 + h_1) G(n)\right| \ll \frac{A}{W^{1/20}} H' X
$$
for any $1$-bounded function $G: \Z \to \C$.

We use a standard ``van der Corput'' argument.  By the Cauchy-Schwarz inequality, it suffices to show that
$$
\sum_{1 \leq h_1 \leq H'} \left|\sum_{1 \leq n \leq X} 1_{\mathcal S}g_1(a_1 n + b_1 + h_1) G(n)\right|^2 \ll \frac{A^2}{W^{1/10}} (H')^2 X^2.
$$
The left-hand side may be rewritten as
$$ \sum_{n,n' \leq X} G(n) \overline{G(n')} \sum_{1 \leq h_1 \leq H'} 1_{\mathcal S}g_1(a_1 n + b_1 + h_1) 1_{\mathcal S}\overline{g_j}(a_1 n' + b_1 + h_1) .$$
By the triangle inequality, it thus suffices to show that
$$ \sum_{n,n' \leq X} \left|\sum_{1 \leq h_1 \leq H'} 1_{\mathcal S}g_{1}(a_1 n + b_1 + h_1) 1_{\mathcal S}\overline{g_1}(a_1 n' + b_1 + h_1)\right|
\ll \frac{A^2}{W^{1/10}} H' X^2.$$
To abbreviate notation we now write $h = h_1$, $g = g_1$, $a = a_1$, $b = b_1$.  By the Cauchy-Schwarz inequality, it suffices to show that
$$ \sum_{n,n' \leq X} \left|\sum_{1 \leq h \leq H'} 1_{\mathcal S}g(a n + b + h) 1_{\mathcal S}\overline{g}(a n' + b + h)\right|^2
\ll \frac{A^4}{W^{1/5}} (H')^2 X^2.$$
Replacing $n,n'$ by $an+b$, $an'+b$ respectively, it suffices to show that
$$ \sum_{n,n'} \left|\sum_{1 \leq h \leq H'} 1_{\mathcal S}g(n + h) 1_{\mathcal S}\overline{g}(n' + h)\right|^2
\ll \frac{A^4}{W^{1/5}} (H')^2 X^2$$
where we have extended $1_{\mathcal S} g$ by zero to the negative integers.
The left-hand side can be rewritten as
$$ \sum_{|h| < H'} (\lfloor H' \rfloor-|h|) \left|\sum_n 1_{\mathcal S} g(n) 1_{\mathcal S} \overline{g}(n+h)\right|^2,$$
and the claim follows from Proposition~\ref{prop:truncEllsc}.
\end{proof}

\appendix

\section{Mean values of complex multiplicative functions in short intervals}\label{complex-app}

In this section we prove a complex variant of results in~\cite{MR} in the case that $f$ is not $p^{it}$ pretentious. In particular we show that the mean value of a $1$-bounded nonpretentious multiplicative function is small for most short intervals:

\begin{theorem}\label{th:MRinCnotInS}
Let $f$ be a $1$-bounded multiplicative function and let $M(f; X)$ be as in~\eqref{eq:MgXdef}. Then, for $X \geq h \geq 10$,
$$
\frac{1}{X}\int_X^{2X} \left|\frac{1}{h} \sum_{\substack{x \leq n \leq x + h}} f(n)\right|^2 dx \ll \exp(- M(f;X)) M(f;X) + \frac{(\log \log h)^2}{(\log h)^2} + \frac{1}{(\log X)^{1/50}}.
$$
\end{theorem}
\begin{remark*}
The factor $\exp(-M(f;X)) M(f;X)$ can be replaced by $\exp(-M(f;X))$, see the remark following Proposition \ref{prop:MeanSquare}
\end{remark*}

Actually as in~\cite{MR} and earlier in this paper, one gets better quantitative results if one first restricts to a subset of $n$ with a typical factorization. Let us first define such subset $\mathcal{S}$ in this setting.

Let $\eta \in (0, 1/6)$, and let $X_0$ be a quantity with $\sqrt{X} \leq X_0 \leq X$.  (The results in \cite{MR} used the choice $X_0=X$, but for technical reasons we will need a more flexible choice of this parameter.)  Consider a sequence of increasing intervals $[P_j, Q_j], j \geq 1$ such that 
\begin{itemize}
\item $Q_1 \leq \exp(\sqrt{\log X_0})$.
\item The intervals are not too far from each other, precisely
\begin{equation}
\label{eq:PjQjnottoofar}
\frac{\log\log Q_{j}}{\log P_{j-1}-1} \leq \frac{\eta}{4j^2}
\end{equation}
for all $j \geq 2$.
\item The intervals are not too close to each other, precisely
\begin{equation}
\label{eq:PjQjnottooclose}
\frac{\eta}{j^2} \log P_{j} \geq 8 \log Q_{j-1} + 16 \log j 
\end{equation}
for all $j \geq 2$.
\end{itemize}
For example, given $0 < \eta < 1/6$, the sequence of intervals $[P_j,Q_j]$ defined in Definition \ref{S-def} can be verified to obey the above estimates if
$$\exp(\sqrt{\log X_0}) \geq Q_1 \geq P_1 \geq (\log Q_1)^{40/\eta}$$
and $P_1$ is sufficiently large.

Let $\mathcal{S}$ be the set of integers $X \leq n \leq 2X$
having at least one prime factor
in each of the intervals $[P_j, Q_j]$ for $j \leq J$, where $J$
is chosen to be the largest index $j$ such that $Q_j \leq \exp((\log X_0)^{1/2})$. 
We will establish the following variant of \cite[Theorem 3]{MR}.

\begin{theorem}\label{th:MRinC}
Let $f$ be a $1$-bounded multiplicative function.
Let $\mathcal{S}$ be as above with $\eta \in (0, 1/6)$. If $[P_1, Q_1] \subset [1, h]$, then for all $X > X(\eta)$ large enough and $h \geq 3$,
$$
\frac{1}{X}\int_X^{2X} \left|\frac{1}{h} \sum_{\substack{x \leq n \leq x + h \\ n \in \mathcal{S}}} f(n)\right|^2 dx \ll \exp(- M(f;X)) M(f;X) + \frac{(\log h)^{1/3}}{P_1^{1/6-\eta}} + \frac{1}{(\log X)^{1/50}}.
$$
\end{theorem}
\begin{remark*}
The factor $\exp(-M(f;X)) M(f;X)$ can be replaced by $\exp(-M(f;X))$, see the remark following Proposition \ref{prop:MeanSquare}
\end{remark*}


The proof of Theorem~\ref{th:MRinC} proceeds as the proof of~\cite[Theorem 3]{MR}. The first step is a Parseval bound
\[
\begin{split}
\frac{1}{X} \int_{X}^{2X} \left|\frac{1}{h} \sum_{\substack{x \leq n \leq x+h \\ n \in \mathcal{S}}} f(n) \right|^2 dx \ll \int_{1}^{1+iX/h_1} \left|F(s)\right|^2 |ds| + \max_{T \geq X/h_1} \frac{X/h_1}{T} \int_{1+iT}^{1+i2T} \left|F(s)\right|^2 |ds|.
\end{split}
\]
This follows exactly in the same way as~\cite[Lemma 14]{MR} but there is no need to split the integral into two parts, and one can just work as for $V(x)$ there. Theorem~\ref{th:MRinC} now follows immediately from the following variant of~\cite[Proposition 1]{MR}.

\begin{proposition}
\label{prop:MeanSquare}
Let $f$ be a $1$-bounded multiplicative
function. Let $\mathcal{S}$ be as above with $\eta \in (0, 1/6)$ and let
\[
F(s) = \sum_{\substack{X \leq n \leq 2X \\ n \in \mathcal{S}}} \frac{f(n)}{n^s}.
\]
Then, for any $T \geq 1$,
\[
\int_{-T}^T \left|F(1+it)\right|^2 dt \ll \left(\frac{T}{X/Q_1} + 1\right) \left(\frac{(\log Q_1)^{1/3}}{P_1^{1/6-\eta}} +  \frac{M(f;X)}{\exp(M(f;X))}  + \frac{1}{(\log X)^{1/50}}\right).
\]
\end{proposition}

\begin{remark*}
In the published version of the paper, the proof of this proposition is incorrect when $M(f;X)$ grows very slowly with $X$. The corrected proof that we provide here gives a slightly stronger result with $\exp(- M(f;X))$ in place of $\exp(- M(f;X))M(f; X)$. We state the result with the weaker factor $\exp(-M(f;X)) M(f; X)$ to remain consistent with the published version of the paper.
\end{remark*}

\begin{proof}
Since the mean value theorem gives the bound $O(T/X + 1)$, we can assume $T \leq X/2$ and $M(f;X) \geq 1$.

Let now $t_1$ be the value of $t$ which attains the minimum in
$$ M(f;X) = \inf_{|t| \leq X} {\mathbb D}(f, n \mapsto n^{it}; X)^2.$$ If $M(f; X) \geq \frac{1}{8} \log \log X$, we write $\mathcal{T}_1 := [-T, T]$ and $\mathcal{T}_0 := \emptyset$ whereas otherwise we write 
\[
\begin{split}
\mathcal{T}_0 &:= \{|t| \leq T \colon |t-t_1| \leq (\log X)^{1/16} \} \\
\mathcal{T}_1 &:= \{|t| \leq T \colon |t-t_1| > (\log X)^{1/16} \}.
\end{split}
\]

Let us first handle $\mathcal{T}_1$. For this we use the following lemma whose proof is in the spirit of works of Granville and Soundararajan (see e.g.~\cite{GSDecay}).
\begin{lemma}
Let $\mathcal{J} \subseteq \{1, \dotsc, J\}$ and $|t| \leq X$. 
\begin{enumerate}[(i)]
\item One has
\begin{equation}
\label{eq:fgdist}
\begin{split}
\mathbb{D}(fg_{\mathcal{J}}, p^{it}; X)^2 \geq \frac{1}{2}\mathbb{D}(f, p^{it}; X)^2
\end{split}
\end{equation}
\item If $M(f; X) \geq \frac{1}{8} \log \log X$ or $|t-t_1| > (\log X)^{1/16}/2$, then
\[
\mathbb{D}(fg_{\mathcal{J}}, p^{it}; X)^2 \geq \left(\frac{1}{6}-\frac{1}{3\pi}-\varepsilon\right) \log \log X
\]
for any $\varepsilon > 0$.
\end{enumerate}
\end{lemma}
\begin{proof}
Let us first show (i). We have
\[
\begin{split}
2\mathbb{D}(fg_{\mathcal{J}}, p^{it}; X)^2 & = \sum_{p \leq X} \frac{1-\Re f(p)g_{\mathcal{J}}(p) p^{-it}}{p} + \sum_{p \leq X} \frac{1-\Re f(p)p^{-it}}{p} + \sum_{p \leq X} \frac{\Re f(p) p^{-it} (1-g_{\mathcal{J}}(p)) }{p} \\
&\geq \sum_{p \leq X} \frac{1-g_{\mathcal{J}}(p)}{p} + \sum_{p \leq X} \frac{1-\Re f(p)p^{-it}}{p} - \sum_{p \leq X} \frac{1-g_{\mathcal{J}}(p)}{p} \geq \mathbb{D}(f, p^{it}; X)^2
\end{split}
\]

Let us now turn to (ii). Notice first that when $M(f; X) \geq \frac{1}{8} \log \log X$, part (i) implies that, whenever $|t| \leq X$, we have $\mathbb{D}(fg_{\mathcal{J}}, p^{it}; X)^2 \geq \frac{1}{16} \log \log X$ which is sufficient. Hence we can concentrate on the case where $M(f; X) < \frac{1}{8} \log \log X$ and $|t-t_1| > (\log X)^{1/16}/2$. 

Writing $Y = \exp((\log X)^{2/3+\varepsilon})$, we have
\begin{equation}
\label{eq:Dfpitlow}
\begin{split}
\mathbb{D}(f, p^{it}; X)^2 &\geq \frac{1}{2} \sum_{p \leq X} \frac{1-\Re f(p)p^{-it}}{p} + \frac{1}{2} \sum_{p \leq X} \frac{1-\Re f(p)p^{-it_1}}{p} \\
&= \sum_{p \leq X} \frac{1-\Re f(p)p^{-i(t+t_1)/2} \cos\left(\frac{(t-t_1) \log p}{2}\right)}{p} \\
&\geq \sum_{Y < p \leq X} \frac{1-\left|\cos\left(\pi \left\Vert \frac{(t-t_1) \log p}{2\pi}\right\Vert\right)\right|}{p},
\end{split}
\end{equation}
where $\Vert x \Vert$ denotes the distance from the nearest integer. When $(\log X)^{1/16}/2 \leq |t-t_1| \leq (\log X)^{20}$, we get as in~\cite[Proof of Lemma 2.3]{GSDecay} by splitting $p$ into short segments $(y, y(1+(\log X)^{-30})]$ that 
\begin{equation}
\label{eq:cossumint}
\begin{split}
\sum_{Y < p \leq X} \frac{1-\left|\cos\left(\pi \left\Vert \frac{(t-t_1) \log p}{2\pi}\right\Vert\right)\right|}{p} &\geq \left(1- \int_0^{1} |\cos (\pi t)| dt \right) \log \frac{\log X}{\log Y} + O(1) \\
&= \left(1-\frac{2}{\pi}\right) \log \frac{\log X}{\log Y} + O(1)
\end{split}
\end{equation}
On the other hand, when $|t-t_1| > (\log X)^{20}$ and $|t| \leq X$
\[
\frac{(t-t_1) \log p}{2\pi}
\]
is equidistributed $\pmod{1}$ by the Erd\H{o}s-Turan inequality and the Vinogradov-Korobov zero free region for $\zeta(s)$ since
\[
e^{2\pi i k\cdot \frac{(t-t_1) \log p}{2\pi}} = p^{i (t-t_1) k}.
\] 
Consequently~\eqref{eq:cossumint} holds also in this case and thus, recalling the definition of $Y$, we obtain from~\eqref{eq:Dfpitlow}
\begin{equation}
\label{eq:fpitdist}
\begin{split}
\mathbb{D}(f, p^{it}; X)^2 &\geq \left(\frac{1}{3}-\frac{2}{3\pi}-\varepsilon\right) \log \log X.
\end{split}
\end{equation}
Now (ii) follows from combining this with part (i).
\end{proof}

The previous lemma implies that for any $t \in \mathcal{T}_1$ we have
\[
\sup_{\substack{t' \\ |t'-t| \leq (\log X)^{1/16}/2}} \mathbb{D}(fg_{\mathcal{J}}, p^{it'}; X)^2 \geq \left(\frac{1}{6}-\frac{1}{3\pi} - \varepsilon\right) \log \log X + O(1).
\]
Hence Hal\'{a}sz's theorem (see e.g.~\cite[Corollary 1 with $T = (\log X)^{1/16}/4$]{GSDecay}) implies that, for any $t \in \mathcal{T}_1$ and any $\varepsilon > 0$
\[
\sum_{\substack{X \leq n \leq 2X}} g_{\mathcal{J}}(n) f(n) n^{-it} \ll \frac{X}{(\log X)^{\frac{1}{6}-\frac{1}{3\pi} - \varepsilon}}.
\]
Proceeding in exactly in the same way as in \cite[Lemma 3]{MR} we obtain the following Lemma. 
\begin{lemma}
\label{le:Halappl}
Let $X \geq Q \geq P \geq 2$. Let $t_1$ be as above, $\varepsilon > 0$ and let $\mathcal{J} \subset \{1, \dotsc, J\}$ and
\[
G(s) = \sum_{\substack{X \leq n \leq 2X}} \frac{g_{\mathcal{J}}(n) f(n)}{n^s} \cdot \frac{1}{\# \{p \in [P, Q] \colon p \mid n\} + 1}.
\]
Then, for any $t \in \mathcal{T}_1$,
\[
|G(1+it)| \ll  \frac{\log Q}{(\log X)^{\frac{1}{6}-\frac{1}{3\pi} - \varepsilon} \log P} + \log X \cdot \exp\left(-\frac{\log X}{3\log Q}\log \frac{\log X}{\log Q} \right).
\]
\end{lemma}
This was the only part in the proof \cite[Proposition 1]{MR} that needed $f$ to be real-valued. In~\cite[Lemma 3]{MR} we had $1/16$ in place of $\rho := \frac{1}{6}-\frac{1}{3\pi} - \varepsilon$ but replacing $1/48$ by $\rho/3 > 1/50$ in the definitions of $P, Q$ and $H$ in the treatment of $\mathcal{U}$ in the proof of~\cite[Proposition 1]{MR}, we still obtain
\[
\int_{\mathcal{T}_1} \left|F(1+it)\right|^2 dt \ll \left(\frac{T}{X/Q_1} + 1\right) \left(\frac{(\log Q_1)^{1/3}}{P_1^{1/6-\eta}} + \frac{1}{(\log X)^{1/50}}\right).
\]

Now we only need to deal with $\mathcal{T}_0$. We shall show that, for $t \in \mathcal{T}_0$ one has 
\begin{equation}
\label{eq:T0claim}
F(1+it) \ll \frac{\exp(-\frac{1}{2} M(f; X))}{1+|t-t_1|} + (\log X)^{-1/16}.
\end{equation}
which immediately implies that 
\[
\int_{\mathcal{T}_0} |F(1+it)|^2 dt \ll \frac{1}{\exp(M(f; X))} + (\log X)^{1/16-2/16}.
\]
Hence Proposition~\ref{prop:MeanSquare} follows once we have shown that~\eqref{eq:T0claim} holds for every $t \in \mathcal{T}_0$. By inclusion-exclusion and partial summation, it suffices to show
\begin{lemma}
\label{le:T0est}
Let $t \in \mathcal{T}_0$. Then
\[
\frac{1}{X} \sum_{\mathcal{J} \subseteq \{1, \dotsc, J\}} (-1)^{\# \mathcal{J}} \sum_{n \leq X} \frac{g_{\mathcal{J}}(n) f(n)}{n^{it}} \ll \frac{\exp(-\frac{1}{2} M(f; X))}{1+|t-t_1|} + (\log X)^{-1/16}.
\]
\end{lemma}

Let us first do the following renormalization:
\begin{lemma}
\label{le:renorm}
Let $t \in \mathcal{T}_0$ and $\mathcal{I} \subseteq \{1, \dotsc, J\}$. Then
\begin{equation}
\label{eq:shift}
\sum_{\substack{n \leq X}} g_{\mathcal{J}}(n) f(n)n^{-it} = \frac{X^{i(t-t_1)}}{1+ i(t-t_1)} \sum_{n \leq X}  g_{\mathcal{J}}(n) f(n)n^{-it_1} + O\left(\frac{X}{(\log X)^{1/10}}\right).
\end{equation}
\end{lemma}

\begin{proof}
Recall that when $t \in \mathcal{T}_0$, we have $|t-t_1| \leq (\log X)^{1/16}$ and $M(f; X) < \frac{1}{8} \log \log X$. We apply \cite[Lemma 7.1]{GSDecay} which gives
\[
\begin{split}
\sum_{\substack{n \leq X}} g_{\mathcal{J}}(n) f(n)n^{-it} &= \frac{X^{i(t_1-t)}}{1+ i(t_1-t)} \sum_{n \leq X}  g_{\mathcal{J}}(n) f(n)n^{-it_1} \\
& \quad + O\left(X \frac{\log\log X}{\log X} \exp \left(\sum_{p \leq X} \frac{|1-f(p)p^{-it_1} g_{\mathcal{J}}(p)|}{p}\right)\right).
\end{split}
\]
Here
\[
\begin{split}
\sum_{p \leq X} \frac{|1-f(p)p^{-it_1} g_{\mathcal{J}}(p)|}{p} & = \sum_{\substack{p \leq X \\ p \in \cup_{j \in \mathcal{J}} (P_j, Q_j]}} \frac{1}{p} + \sum_{\substack{p \leq X \\ p \not \in \cup_{j \in \mathcal{J}} (P_j, Q_j]}} \frac{|1-f(p)p^{-it_1}|}{p}
\end{split}
\]
As in~\cite[Proof of Corollary 3]{GSDecay}, note that for $z$ in the unit disc, $|1-z| = (1 + |z|^2 - 2 \Re z)^{1/2} \leq \left(2-2\Re z\right)^{1/2}$. Applying this and Cauchy-Schwarz, we see that
\[
\begin{split}
&\sum_{\substack{p \leq X \\ p \not \in \cup_{j \in \mathcal{J}} (P_j, Q_j]}} \frac{|1-f(p)p^{-it_1}|}{p} \leq \sum_{\substack{p \leq X \\ p \not \in \cup_{j \in \mathcal{J}} (P_j, Q_j]}} \frac{\sqrt{2-2 \Re f(p)p^{-it_1}}}{p} \\
&\leq \sqrt{\sum_{\substack{p \leq X \\ p \not \in \cup_{j \in \mathcal{J}} (P_j, Q_j]}} \frac{2}{p}} \cdot \sqrt{ \sum_{\substack{p \leq X}} \frac{1-\Re f(p)p^{-it_1}}{p}} = \sqrt{\sum_{\substack{p \leq X \\ p \not \in \cup_{j \in \mathcal{J}} (P_j, Q_j]}} \frac{1}{p}} \cdot \sqrt{2M(f; X)}.
\end{split}
\]
Define $\beta$ by 
\[
\beta \log \log X = \sum_{\substack{p \leq X \\ p \in \cup_{j \in \mathcal{J}} (P_j, Q_j]}} \frac{1}{p}.
\]
Since $Q_J \leq \exp((\log X)^{1/2})$, necessarily $\beta \in [0, 1/2 + O(1/\log \log X)]$. Recalling $M(f; X) < \frac{1}{8} \log \log X$, we obtain
\[
\sum_{p \leq X} \frac{|1-f(p)p^{-it_1} g_{\mathcal{J}}(p)|}{p} \leq \beta \log \log X+ \frac{\sqrt{1-\beta}}{2} \log \log X + O(\sqrt{\log \log X}). 
\]
It is easy to see that the right hand side is increasing in $\beta$ in our range, so it is maximized when $\beta = 1/2 + O(1/\log \log X)$ in which case one gets a bound that is $\leq \frac{7}{8} \log \log X$ and the claim follows.
\end{proof}

Notice that $2^J (\log X)^{-1/10} \ll (\log X)^{-1/10+o(1)}$. Hence, thanks to Lemma~\ref{le:renorm}, Lemma~\ref{le:T0est} follows  once we have shown that
\begin{equation}
\label{eq:normclaim}
U := \frac{1}{X} \sum_{\mathcal{J} \subseteq \{1, \dotsc, J\}} (-1)^{\# \mathcal{J}} \sum_{n \leq X} g_{\mathcal{J}}(n) f(n) n^{-it_1}  \ll \exp\left(-\frac{1}{2} M(f; X)\right) + (\log X)^{-1/16}.
\end{equation}

We use the method of the proof of Hal\'{a}sz's theorem from~\cite{GHS}. Take $T_0 := (\log X)^2$ and $y:= Q_J$. 
Define $s_{\mathcal{J}}$ and $\ell$ to be the multiplicative functions with 
\[
s_{\mathcal{J}}(p^k) = f(p^k) g_{\mathcal{J}}(p^k) \mathbf{1}_{p \leq y} \quad \text{and} \quad \ell(p^k) = f(p^k) \mathbf{1}_{p > y}
\] 
and set
\[
\mathcal{S}_{\mathcal{J}}(s) = \sum_{n \geq 1} \frac{s_{\mathcal{J}}(n)}{n^s} \quad \text{and} \quad \mathcal{L}(s) = \sum_{n \geq 1} \frac{\ell(n)}{n^s}.
\]
Furthermore, define $\Lambda_\ell(n)$ through
\[
-\frac{\mathcal{L}'(s)}{\mathcal{L}(s)} = \sum_{n = 2}^\infty \frac{\Lambda_\ell(n)}{n^s}.
\]
Write $\eta = 1/\log y$ and $c_0 = 1+1/\log X$. Now we apply~\cite[Proposition 2.1]{GHS}; strictly speaking we do not necessarily have $|\Lambda_f(n)| \leq \Lambda(n)$ as required when using arguments of~\cite{GHS} with $\kappa = 1$, but this inequality does hold for square-free $n$ and it is easy to see that the relevant parts of~\cite{GHS} with $\kappa = 1$ work for all $1$-bounded multiplicative functions. We get
\begin{equation}
\label{eq:HGSProp2.1}
\begin{split}
&\sum_{n \leq X} f(n) g_{\mathcal{J}}(n)n^{-it_1} = \int_0^\eta \int_0^\eta \frac{1}{\pi i} \int_{c_0-iT_0}^{c_0+iT_0} \mathcal{S}_{\mathcal{J}} (s-\alpha - \beta + it_1) \mathcal{L}(s+\beta+it_1) \\
& \qquad \cdot \sum_{\substack{y < m < x/y \\ y < n < x/y}} \frac{\Lambda_{\ell}(m)}{m^{s+it_1-\beta}} \frac{\Lambda_{\ell}(n)}{n^{s+it_1+\beta}} \frac{X^{s - \alpha-\beta}}{s-\alpha - \beta} ds d\beta d\alpha + O\left(\frac{X}{\log X} \log y \right).
\end{split}
\end{equation}
Write $s = c_0 + it_0$ and $c_0+it_0-\alpha-\beta + it_1 = \sigma + it$ with $\sigma, t \in \mathbb{R}$. Note that on the right hand side only $\mathcal{S}_{\mathcal{J}}(\sigma + it)$ depends on $\mathcal{J}$ and furthermore
\begin{equation}
\label{eq:SJtreat}
\begin{split}
&\sum_{\mathcal{J} \subseteq \{1, \dotsc, J\}} (-1)^{\# \mathcal{J}} \mathcal{S}_{\mathcal{J}}(\sigma + it) = \sum_{\substack{n \\ p \mid n \implies p \leq y}} \frac{f(n)}{n^{\sigma + it}} \prod_{j = 1}^J (1-g_{j}(n))\\
& = \prod_{\substack{p \leq y \\ p \not \in \cup_{j=1}^J (P_j, Q_j]}} \left(1+\frac{f(p)}{p^{\sigma +it}} + \dotsb\right) \prod_{j = 1}^J \left(\prod_{p \in (P_j, Q_j]} \left(1+\frac{f(p)}{p^{\sigma +it}} + \dotsb\right) - 1\right) \\
&\ll \left|\prod_{p \leq y} \left(1+\frac{f(p)}{p^{\sigma+it}}\right)^{1/2} \left(1+\frac{1}{p^\sigma}\right)^{1/2}\right| \\
& \cdot \prod_{j= 1}^J \left|\exp\left(\frac{1}{2}\sum_{P_j < p \leq Q_j} \frac{f(p)}{p^{\sigma + it}}\right) - \exp\left(-\frac{1}{2}\sum_{P_j < p \leq Q_j} \frac{f(p)}{p^{\sigma +it}}\right)\right| \exp\left(-\frac{1}{2} \sum_{P_j < p \leq Q_j} \frac{1}{p^\sigma}\right).
\end{split}
\end{equation}
Next we shall show that the expression on the last line is always $\leq 1$. For $j = 1, \dotsc, J$, define real numbers
\[
\alpha_j := \Re \sum_{P_j < p \leq Q_j} \frac{f(p)}{p^{\sigma +it}} \quad \text{and} \quad \theta_j := \Im \sum_{P_j < p \leq Q_j} \frac{f(p)}{p^{\sigma+it}},
\]
so that we wish to show that, for every $j = 1, \dotsc, J$,
\[
\left|\exp\left(\frac{1}{2}(\alpha_j + i\theta_j)\right) - \exp\left(-\frac{1}{2}(\alpha_j + i\theta_j)\right)\right| \leq \exp\left(\frac{1}{2}\sum_{P_j < p \leq Q_j} \frac{1}{p^\sigma}\right).
\]
By definition and triangle inequality
\[
\alpha_j^2 + \theta_j^2 = \left|\sum_{P_j < p \leq Q_j} \frac{f(p)}{p^{\sigma+it}}\right|^2 \leq \left(\sum_{P_j < p \leq Q_j} \frac{1}{p^\sigma}\right)^2
\]
and thus it suffices to show that
\begin{equation}
\label{eq:leq1}
\left|\exp\left(\frac{1}{2}(\alpha_j + i\theta_j)\right) - \exp\left(-\frac{1}{2}(\alpha_j + i\theta_j)\right)\right| \leq \exp\left(\frac{1}{2}\sqrt{\alpha_j^2 + \theta_j^2}\right).
\end{equation}
Squaring both sides and applying the law of cosines, this reduces to
\begin{lemma}
Let $\alpha, \theta \in \mathbb{R}$. Then
\[
e^\alpha + e^{-\alpha} -  2 \cos \theta \leq e^{\sqrt{\alpha^2 + \theta^2}}.
\]
\end{lemma}
\begin{proof}
By symmetry we can assume that $\alpha, \theta > 0$. The function $x \mapsto e^{\sqrt{x}}$ has derivative $\frac{1}{2} x^{-1/2} e^{\sqrt{x}}$. Differentiating again, we see that this derivative is minimised at $x=1$ with value $e/2$, so by the mean value theorem
\[
e^{\sqrt{\alpha^2 + \theta^2}} \geq e^\alpha + \frac{e}{2} \theta^2
\]
On the other hand
\[
\cos \theta = 1 - 2 \sin^2 \frac{\theta}{2} \geq 1 - \frac{1}{2} \theta^2.
\]
So it suffices to show that
\[
e^\alpha + e^{-\alpha} -2 + \theta^2 \leq e^\alpha + \frac{e}{2} \theta^2 \iff 2 - e^{-\alpha} +\theta^2\left(\frac{e}{2}-1\right) \geq 0.
\]
But this follows immediately since $2 \geq e^{-\alpha}$ and $e/2 \geq 1$.
\end{proof}
Hence we indeed got that the expression on the last line of~\eqref{eq:SJtreat} is always $\leq 1$. Using also that $\sigma = c_0 -\alpha-\beta$ and $\alpha, \beta < 1/\log y$, we see that
\begin{equation}
\label{eq:Sestimate}
\begin{split} 
\sum_{\mathcal{J} \subseteq \{1, \dotsc, J\}} (-1)^{\# \mathcal{J}} \mathcal{S}_{\mathcal{J}}(\sigma + it) &\ll |\mathcal{S}_{\emptyset}(\sigma + it)|^{1/2} \prod_{\substack{p \leq y}} \left(1+\frac{1}{p^\sigma}\right)^{1/2} \\
&\ll |\mathcal{S}_{\emptyset}(c_0 + \beta + it_0 + it_1)|^{1/2} \prod_{\substack{p \leq y}} \left(1+\frac{1}{p^{c_0 + \beta}}\right)^{1/2}.
\end{split}
\end{equation}
Furthermore
\begin{equation}
\label{eq:Lestimate}
\begin{split}
&\mathcal{L}(c_0+\beta+it_0+it_1) \ll \left|\mathcal{L}(c_0+\beta+it_0+it_1)\right|^{1/2} \prod_{\substack{p > y}} \left(1+\frac{1}{p^{c_0 + \beta}}\right)^{1/2}.
\end{split}
\end{equation}
Write 
\[
F_0(s) = \mathcal{S}_\emptyset(s) \mathcal{L}(s) = \sum_n \frac{f(n)}{n^s}.
\]
Recall that we aim to prove~\eqref{eq:normclaim}. Plugging~\eqref{eq:Sestimate} and~\eqref{eq:Lestimate} into~\eqref{eq:HGSProp2.1} and rearranging, we see that
\[
\begin{split}
U &\ll \int_0^\eta \int_0^\eta X^{1-\alpha-\beta} \max_{|t_0| \leq T_0} \frac{|F_0(c_0+\beta + it_0+ it_1)|^{1/2} \zeta(c_0 + \beta)^{1/2}}{|c_0 + \beta + it_0|} \\
&\qquad \cdot \int_{c_0-iT_0}^{c_0+iT_0} \left|\sum_{\substack{y < m < x/y \\ y < n < x/y}} \frac{\Lambda_{\ell}(m)}{m^{s+it_1-\beta}} \frac{\Lambda_{\ell}(n)}{n^{s+it_1+\beta}}\right| |ds| d\beta d\alpha + \frac{X}{(\log X)^{1/2}}.
\end{split}
\]
As in~\cite[Proof of Theorem 1.1]{GHS}, we get from this that
\[
U \ll \frac{X}{\log X} \int_{1/\log X}^{2/\log y} \max_{|t_0| \leq T_0} \frac{|F_0(1+\sigma' + it_0+ it_1)|^{1/2} \zeta(1+\sigma')^{1/2}}{|1+\sigma' + it_0 + it_1|} \frac{d\sigma'}{\sigma'} + \frac{X}{(\log X)^{1/2}}.
\]
Using the maximum modulus principle as in~\cite[Proof of Corollary 1.2]{GHS}, we see that
\[
\begin{split}
\max_{\substack{\tfrac{1}{\log X} \leq \sigma' \leq 2 \\ |t'-t_1| \leq T_0}} \frac{|F_0(1+\sigma' + it')|^{1/2} }{|1+\sigma' + it'|} \ll \max_{|t'-t_1| \leq T_0} |F_0(1+\tfrac{1}{\log X} + it')|^{1/2} + 1 \ll \frac{(\log X)^{1/2}}{\exp(M(f; X)/2)}  + 1.
\end{split}
\]
Since $\zeta(1+\sigma') \ll 1/\sigma'$, we obtain
\[
\begin{split}
U &\ll \frac{X}{\log X} \left(\frac{(\log X)^{1/2}}{\exp(M(f; X)/2)}  + 1\right) \int_{1/\log X}^{2/\log y} \frac{d\sigma'}{\sigma'^{3/2}} + \frac{X}{(\log X)^{1/2}} \\
&\ll \frac{X}{\exp(M(f; X)/2)} + \frac{X}{(\log X)^{1/2}}.
\end{split}
\]
Hence~\eqref{eq:normclaim} holds so that the proof of Lemma~\ref{le:T0est} is finished. As described below, this implies Proposition~\ref{prop:MeanSquare} and thus also Theorem~\ref{th:MRinC}.
\end{proof}

\begin{proof}[Proof of Theorem~\ref{th:MRinCnotInS}]
Let $\eta = 1/12$, $P_1 = (\log h)^{480}, Q_1 = h$, let $P_j$ and $Q_j$ for $j \geq 2$ be as in Definition~\ref{S-def}, and let $\mathcal{S}$ be as above. Then
\[
\frac{1}{X}\int_X^{2X} \left|\frac{1}{h} \sum_{\substack{x \leq n \leq x + h}} f(n)\right|^2 dx \leq \frac{1}{X}\int_X^{2X} \left|\frac{1}{h} \sum_{\substack{x \leq n \leq x + h \\ n \in \mathcal{S}}} f(n)\right|^2 dx + \frac{1}{X}\int_X^{2X} \left|\frac{1}{h} \sum_{\substack{x \leq n \leq x + h \\ n \not \in \mathcal{S}}} 1 \right|^2 dx.
\]
The contribution from the first integral is acceptable by Theorem~\ref{th:MRinC}.
We rewrite the second integrand as
\[
\begin{split}
\left|\frac{1}{h} \sum_{\substack{x \leq n \leq x + h \\ n \not \in \mathcal{S}}} 1\right| &= \left|1+O(1/h) - \frac{1}{h} \sum_{\substack{x \leq n \leq x + h \\ n \in \mathcal{S}}} 1 \right| \\
&\leq \left|\frac{1}{X} \sum_{\substack{X \leq n \leq 2X \\ n \in \mathcal{S}}} 1 - \frac{1}{h} \sum_{\substack{x \leq n \leq x + h \\ n \in \mathcal{S}}} 1\right| + \left|\frac{1}{X} \sum_{\substack{X \leq n \leq 2X \\ n \not \in \mathcal{S}}} 1\right| + O(1/h),
\end{split}
\]
and the claim follows from~\cite[Theorem 3 with $f=1$]{MR} and Lemma~\ref{excep}.
\end{proof}

\section{Counterexample to the uncorrected Elliott conjecture}\label{counter}

In this appendix we present a counterexample to Conjecture \ref{elliott-conj}.  More precisely:

\begin{theorem}[Counterexample]  There exists a $1$-bounded multiplicative function $g: \N \to \C$ such that
\begin{equation}\label{gcp}
 \sum_p \frac{1 - \operatorname{Re}( g(p) \overline{\chi(p)} p^{-it} )}{p} = \infty
\end{equation}
for all Dirichlet characters $\chi$ and $t \in \R$ (i.e., one has $M(g;\infty,\infty)=\infty$), but such that
\begin{equation}\label{lan}
\left|\sum_{n \leq t_m} g(n) \overline{g(n+1)}\right| \gg t_m
\end{equation}
for all sufficiently large $m$, and some sequence $t_m$ going to infinity.
\end{theorem}

\begin{proof}  For each prime $p$, we choose $g(p)$ from the unit circle $S^1 := \{ z: |z|=1\}$ by the following iterative procedure involving a sequence $t_1 < t_2 < t_3 < \dots$.

\begin{enumerate}
\item Initialize $t_1 := 100$ and $m := 1$, and set $g(p):=1$ for all $p \leq t_1$.
\item Now suppose recursively that $g(p)$ has been chosen for all $p \leq t_{m}$.  As the quantities $\log p$ are linearly independent over the integers, the (continuous) sequence $t \mapsto (t \log p \hbox{ mod } 1)_{p \leq t_{m}}$ is equidistributed in the torus $\prod_{p \leq t_{m}} {\mathbb T}$; equivalently, the sequence $t \mapsto (p^{it})_{p \leq t_{m}}$ is equidistributed in the torus $\prod_{p \leq t_{m}} S^1$.  Thus one can find a quantity $s_{m+1} > \exp(t_m)$ such that
\begin{equation}\label{gap}
p^{is_{m+1}} = g(p) \left(1 + O\left( \frac{1}{t_m^2} \right)\right)
\end{equation}
for all $p \leq t_m$.
\item Set $t_{m+1} := s_{m+1}^2$, and then set 
\begin{equation}\label{gap-2}
g(p) := p^{is_{m+1}}
\end{equation}
 for all $t_m < p \leq t_{m+1}$.  Now increment $m$ to $m+1$ and return to step 2.
\end{enumerate}

Clearly the $t_m$ go to infinity, so $g(p)$ is defined for all primes $p$.
We then define 
\begin{equation}\label{gan}
 g(n) := \mu(n)^2 \prod_{p|n} g(p),
\end{equation}
which is clearly a $1$-bounded multiplicative function.

Suppose that $n \leq t_{m+1}$ is squarefree. Then $n$ is the product of distinct primes less than or equal to $t_{m+1}$, including at most $t_m$ primes less than or equal to $t_m$.  From \eqref{gan} we then have
\begin{align*}
g(n) &= n^{is_{m+1}} \left(1 + O\left( \frac{1}{t_m^2} \right)\right)^{O( t_m )} \\
&= n^{is_{m+1}} + O\left( \frac{1}{t_m} \right).
\end{align*}
If $n$ is not squarefree, then $g(n)$ of course vanishes.  We thus have, for $t_{m+1}^{3/4} \leq n \leq t_{m+1}-1$,
\begin{align*}
 g(n) \overline{g(n+1)} &= \mu^2(n) \mu^2(n+1) \left( \frac{n+1}{n} \right)^{i s_{m+1}} + O\left( \frac{1}{t_m} \right) \\
&= \mu^2(n) \mu^2(n+1) + O\left( \frac{s_{m+1}}{t_{m+1}^{3/4}} \right) + O\left( \frac{1}{t_m} \right) \\
&= \mu^2(n) \mu^2(n+1) + O\left( \frac{1}{t_m} \right),
\end{align*}
and the claim \eqref{lan} then easily follows since the sequence $\mu^2(n) \mu^2(n+1)$ has positive mean value.

Now we prove \eqref{gcp}. From \eqref{gap-2}, we have
\[
\begin{split}
\sum_p \frac{1 - \operatorname{Re}( g(p) \overline{\chi(p)} p^{-it} )}{p} &\geq \sum_{t_m < p \leq t_{m+1}} \frac{1 - \operatorname{Re}( \overline{\chi(p)} p^{i(s_{m+1}-t)} )}{p} \\ &\geq \sum_{\exp((\log t_{m+1})^{5/6}) < p \leq t_{m+1}} \frac{1 - \operatorname{Re}( \overline{\chi(p)} p^{i(s_{m+1}-t)} )}{p}
\end{split}
\]
since $\exp((\log t_{m+1})^{5/6}) \geq \exp((2t_m)^{5/6}) \geq t_m$. Hence we see as in \eqref{eq:MestLiou} that the right-hand side goes to infinity as $m \to \infty$ for any fixed $\chi,t$, and the claim follows.
\end{proof}

It is easy to see that the function $g$ constructed in the above counterexample violates \eqref{geo}, and so is not a counterexample to the corrected form of Conjecture \ref{elliott-conj}.  It is also not difficult to modify the above counterexample so that the function $g$ is completely multiplicative instead of multiplicative, using the fact that most numbers up to $t_{m+1}$ have fewer than $t_m$ prime factors less than $t_m$ (counting multiplicity); we leave the details to the interested reader.

\section{An argument of Granville and Soundararajan}\label{gs-sec}

In this appendix we show the equivalence of the hypotheses \eqref{mgi} and \eqref{geo} for Elliott's conjecture in the case that the multiplicative function $g_{j_0}$ is real.  The key lemma is the following estimate, essentially due to Granville and Soundararajan.

\begin{lemma}
\label{le:distest}
Let $f \colon \mathbb{N} \to [-1, 1]$ be a multiplicative function, let $x \geq 100$, and let $\chi$ be a fixed Dirichlet character. For $1 \leq |\alpha| \leq x$, one has
\begin{equation}
\label{eq:firsttrineq} 
{\mathbb D}( f, n \mapsto \chi(n) n^{i\alpha}; x ) \geq \frac{1}{4}\sqrt{\log \log x} + O_\chi(1).
\end{equation}
When $\chi^2$ is non-principal, this holds for all $|\alpha| \leq x$.

If $\chi^2$ is principal (i.e., $\chi$ is a quadratic character), then, for $|\alpha| \leq 1$, one has
\begin{equation}
\label{eq:sectrineq} 
{\mathbb D}( f, n \mapsto \chi(n) n^{i\alpha}; x ) \geq \frac{1}{3} \mathbb{D}(f, \chi ;x) + O(1). 
\end{equation}
\end{lemma}

\begin{proof}
To establish~\eqref{eq:firsttrineq}, we notice that, by conjugation symmetry and the triangle inequality,
\begin{align*}
{\mathbb D}( f, n \mapsto \chi(n) n^{i\alpha}; x ) &= \frac{1}{2}( {\mathbb D}( f, n \mapsto \chi(n) n^{i\alpha}; x ) + {\mathbb D}( f, n \mapsto \overline{\chi(n)} n^{-i\alpha}; x ))\\
&\geq \frac{1}{2}{\mathbb D}( n \mapsto \overline{\chi(n)} n^{-i\alpha}, n \mapsto \chi(n) n^{i\alpha}; x ) \\
&= \frac{1}{2} \left( \sum_{p \leq x} \frac{1 - \operatorname{Re} \chi^2(p) p^{2i\alpha}}{p} \right)^{1/2}
\end{align*}
which implies the claim for $|\alpha| \geq 1$ or for non-principal $\chi^2$ by the zero-free (and pole-free) region for Dirichlet $L$-functions (see~(\ref{eq:MestLiou}) for a related argument).

To establish~\eqref{eq:sectrineq}, notice first that since $\chi^2$ is principal, $\chi$ is real-valued which implies together with the triangle inequality
\[
\mathbb{D}(f, n \mapsto \chi(n) n^{i\alpha}; x) = \mathbb{D}(f \chi, n \mapsto n^{i\alpha}; x) \geq \mathbb{D}(1, f \chi; x) - \mathbb{D}(1, n \mapsto n^{i\alpha}; x).
\]
Now $\mathbb{D}(1, n \mapsto  n^{i\alpha}; x) = \mathbb{D}(1, n \mapsto  n^{2i\alpha}; x) + O(1)$ for $|\alpha| \leq 1$, since $\mathbb{D}(1, n \mapsto  n^{i\alpha}; x)^2 = \log (1 + |\alpha| \log x) + O(1)$ from the prime number theorem, so that the claim follows unless $\mathbb{D}(1, n \mapsto  n^{2i\alpha}; x) \geq \tfrac{2}{3} \mathbb{D}(1, f\chi; x)$. But in the latter case, the triangle inequality gives
\begin{align*}
\frac{2}{3} \mathbb{D} (f,\chi; x) 
&= \frac{2}{3} \mathbb{D} (1,f\chi; x) \\
& \leq \mathbb{D}(1, n \mapsto n^{2i\alpha}; x) \\
&= \mathbb{D}(n \mapsto n^{-i\alpha}, n \mapsto n^{i\alpha}; x) \\ & \leq \mathbb{D}(f\chi, n \mapsto n^{-i\alpha}; x)
+ \mathbb{D}(f\chi; n \mapsto n^{i\alpha}; x) \\
&= 2 \mathbb{D}(f, n \mapsto \chi(n) n^{i\alpha}; x),
\end{align*}
and the claim \eqref{eq:sectrineq} follows.
\end{proof}

From this lemma, we see that when $g_{j_0}$ is a real $1$-bounded multiplicative function, then for given $Q$, the condition \eqref{geo} is equivalent to
$$ {\mathbb D}(g_{j_0}, \chi; X) \to \infty $$ when $X \to \infty$ for all quadratic characters $\chi$ of modulus at most $Q$. But this follows from \eqref{mgi}.  The converse implication is trivial.


\bibliographystyle{plain}
\bibliography{Chowla}

\end{document}